\documentclass[letterpaper,10pt,twoside,journal]{IEEEtran}

\usepackage{epsfig}
\usepackage{bm}
\usepackage{graphics}
\usepackage{amsmath,amssymb,amsfonts,amsbsy,amsthm}
\usepackage{footnote}
\usepackage{multicol}
\usepackage{multirow}
\usepackage{url}
\usepackage{indentfirst}
\usepackage{longtable}
\usepackage{algorithm2e}
\usepackage{array}
\usepackage{mdwmath}
\usepackage{mdwtab}
\usepackage{rotating}
\usepackage{color}

\newcommand {\beq}{\begin{equation}}
\newcommand {\eeq}{\end{equation}}
\newcommand {\beqn}{\begin{eqnarray}}
\newcommand {\eeqn}{\end{eqnarray}}
\newcommand {\bsp}{\begin{split}}
\newcommand {\esp}{\end{split}}
\newcommand {\bsmtx}{\begin{smallmatrix}}
\newcommand {\esmtx}{\end{smallmatrix}}
\newcommand {\bpmtx}{\begin{pmatrix}}
\newcommand {\epmtx}{\end{pmatrix}}
\newcommand {\bbmtx}{\begin{bmatrix}}
\newcommand {\ebmtx}{\end{bmatrix}}
\newcommand {\bBmtx}{\begin{Bmatrix}}
\newcommand {\eBmtx}{\end{Bmatrix}}

\def\m{\mathbf}
\def\mc{\mathcal}
\def\mx{\mathbf{x}}
\def\mr{\mathbf{r}}

\def\mA{\mathbf{A}}

\def\my{\mathbf{y}}
\def\mz{\mathbf{z}}
\def\mX{\mathbf{X}}
\def\mY{\mathbf{Y}}
\def\mU{\mathbf{U}}
\def\mV{\mathbf{V}}
\def\mZ{\mathbf{Z}}
\def\mS{\bm{\Sigma}}

\def\ms{\bm{\sigma}}

\def\mP{\bm{\Psi}}
\def\mW{\bm{\Omega}}
\def\ml{\bm{\lambda}}

\def\l{\left}
\def\r{\right}

\def\lg{\lambda}
\def\Lg{\Lambda}
\def\w{\omega}
\def\W{\Omega}

\def\a{\alpha}
\def\b{\beta}
\def\P{\Psi}
\def\p{\psi}
\def\t{\theta}

\def\suml{\sum\limits}

\newcommand{\br}[1]{\l(#1\r)}
\newcommand{\pr}[1]{\l[#1\r]}
\newcommand{\cbr}[1]{\l\{#1\r\}}
\newcommand{\norm}[2]{\ensuremath{\l\|#1\r\|_{#2}}}
\newcommand{\snorm}[2]{\ensuremath{\l\|#1\r\|_{\mc{S}_{#2}}}}

\newcommand{\abs}[1]{\ensuremath{\l|#1\r|}}

\newcommand{\grad}{\ensuremath{\nabla}}

\DeclareMathOperator*{\argmin}{arg\,min}
\DeclareMathOperator*{\argmax}{arg\,max}
\newcommand{\ip}[2]{\ensuremath{\l<#1\,,\,#2\r>}} 
\def\transpose{T}
\def\hermitian{H}
\newcommand{\Dxx}{\ensuremath{\Delta_{r_1r_1}}}
\newcommand{\Dxy}{\ensuremath{\Delta_{r_1r_2}}}
\newcommand{\Dyy}{\ensuremath{\Delta_{r_2r_2}}}
\def\mb{\mbox}
\newcommand{\sgn}[1]{\ensuremath{\mbox{sgn}\br{#1}}}
\newcommand{\diag}[1]{\ensuremath{\mbox{diag}\br{#1}}}
\newtheorem{lemma}{Lemma}
\newtheorem{proposition}{Proposition}

\newtheorem{thm}{Theorem}
\newtheorem{definition}{Definition}
\def\S{\mathbb{S}}
\def\R{\mathbb{R}}
\def\U{\mathbb{U}}
\def\C{\mathbb{C}}

\def\D{\mathbb{D}}

\newcommand{\dH}{\ensuremath{\mbox{\boldmath{$\mc{H}$}}}} 
\newcommand{\tr}[1]{\ensuremath{\operatorname{tr}\br{#1}}} 

\begin{document}

\urldef{\slEmail}\url{stamatis.lefkimmiatis@epfl.ch}
\urldef{\muEmail}\url{michael.unser@epfl.ch}
\urldef{\jwEmail}\url{john.ward@epfl.ch}

\urldef{\ourEmails}\url{{stamatis.lefkimmiatis,john.ward,michael.unser}@epfl.ch}
\urldef{\ourWebPage}\url{http://bigwww.epfl.ch/}

\title{Hessian Schatten-Norm Regularization \\for  Linear Inverse Problems}

\author{Stamatios~Lefkimmiatis,~\IEEEmembership{Member,~IEEE},
     John~Paul~Ward,  and~Michael~Unser,~\IEEEmembership{Fellow,~IEEE}%
   \thanks{This work was supported (in part) by the Hasler Foundation and the Indo-Swiss Joint Research Program.}
 \thanks{ Copyright (c) 2012 IEEE. Personal use of this material is permitted. However, permission to use this material for any other purposes must be obtained from the IEEE by sending a request to pubs-permissions@ieee.org.}
  \thanks{The authors are with the Biomedical Imaging Group (BIG), \'{E}cole polytechnique f\'{e}d\'{e}rale de Lausanne (EPFL), CH-1015 Lausanne, Switzerland (email: \slEmail;\jwEmail;\muEmail).}%
  \thanks{Digital Object Identifier 10.1109/TIP.2013.2237919}}

\markboth{To appear in {IEEE} Transactions on Image Processing}{Lefkimmiatis \MakeLowercase{\textit{et al.}}: Hessian Schatten-Norm Regularization for Linear Inverse Problems}

\maketitle

\begin{abstract}
We introduce a novel family of invariant, convex, and non-quadratic functionals that we employ to derive regularized solutions of ill-posed linear inverse imaging problems. The proposed regularizers involve the Schatten norms of the Hessian matrix, computed at every pixel of the image. They can be viewed as second-order extensions of the popular total-variation (TV) semi-norm since they satisfy the same invariance properties. Meanwhile, by taking advantage of second-order derivatives, they avoid the staircase effect, a common artifact of TV-based reconstructions, and perform well for a wide range of applications.  To solve the corresponding optimization problems, we propose an algorithm that is based on a primal-dual formulation. A fundamental ingredient of this algorithm is the projection of matrices onto Schatten norm balls of arbitrary radius. This operation is performed efficiently based on a direct link we provide between vector projections onto $\ell_q$ norm balls and matrix projections onto Schatten norm balls. Finally, we demonstrate the effectiveness of the proposed methods through experimental results on several inverse imaging problems with real and simulated data.
\end{abstract}

\ifCLASSOPTIONjournal
\fi

\section{Introduction}
\label{sec:intro}
\IEEEPARstart{L}{inear} inverse problems arise in a host of imaging applications, ranging from microscopy and medical imaging to remote sensing and astronomical imaging~\cite{Bertero1998}.  The task is to reconstruct the underlying image from a series of degraded measurements. These problems are often formulated within a variational framework, where image reconstruction can be cast as the minimization of an energy functional subject to some penalty. The role of the penalty is significant, since it imposes certain constraints on the solution and considerably affects the quality of the reconstruction.

The importance of choosing an appropriate penalty has initiated the development of regularization functionals that can effectively model certain properties of natural images. A popular regularization criterion is the total-variation (TV) semi-norm~\cite{Rudin1992} which has been successfully applied to several imaging problems such as  image denoising,  restoration~\cite{Bioucas2006c,Beck2009b}, inpainting~\cite{Chan2002}, zooming~\cite{Chambolle2004}, and MRI reconstruction~\cite{Lustig2007}. TV owes its success to its ability to preserve the edges of the underlying image well.  Its downside, however, is that it introduces blocking artifacts (a.k.a. \emph{staircase effect})~\cite{Chan2000}. The reason is that TV favors vanishing first-order derivatives. Thus, it tends to result in piecewise-constant solutions even when the underlying images are not necessarily piecewise constant. This tendency  is responsible for oversharpening the contrast along image contours and can be a serious drawback  in many applications.

A common workaround to prevent the oversharpening of regions with smooth intensity transitions is to replace TV by functionals that involve higher-order differential operators, because higher-order derivatives can potentially restore a wider class of images. Often, moving from piecewise-constant to piecewise-linear reconstructions offers a satisfactory improvement in the fitting of smooth intensity changes, so that most of the published functionals involve second-order differentials. Such regularizers have been considered, mostly for image denoising, either combined with TV~\cite{Chan2000}--\cite{Setzer2011}
or in a standalone way~\cite{You2000}--\cite{Bredies2010}.
These recent  advances motivate us to investigate a class of regularizers that depend on matrix norms of the Hessian. These regularizers enjoy most of the favorable properties of TV; namely, convexity, contrast, rotation, translation, and scale invariance (up to a multiplicative constant), while they avoid the staircase effect by not penalizing first-order polynomials.

The key contributions of this work are as follows:
\begin{enumerate}
\item The identification of a novel family of invariant functionals that involve Schatten norms of the Hessian matrix, computed at every pixel of the image. These are used in a variational framework to derive regularized solutions of ill-posed linear inverse imaging problems. Our functionals capture curvature information related to the image intensity and lead to reconstructions that avoid the staircase effect. 
\item A general first-order algorithm for solving the resulting constrained optimization problems under any choice of Schatten norm. The proposed algorithm relies on our derivation of a primal-dual formulation and copes with the non-smooth nature and the high dimensionality of the problem.
\item A direct link between matrix projections onto Schatten norm balls and vector projections onto $\ell_q$ norm balls. This link enables us to design an efficient method for performing matrix projections. Although it is  a fundamental component of our optimization algorithm, our result is not specific to the Hessian and can potentially have a wider applicability.
 \end{enumerate}

The rest of the paper is organized as follows: In Section~\ref{sec:HOR}, we discuss regularization functionals that are commonly used in imaging problems. Then, by focusing on invariance principles we derive our novel family of non-quadratic second-order functionals. In Section~\ref{sec:image_reconstruction}, we present the discrete formulation of the problem and we describe the proposed optimization algorithm. In Section~\ref{sec:results}, we assess the performance of our approach for several linear inverse imaging problems with experiments on standard test images and real biomedical images. We conclude our work in Section~\ref{sec:conclusions}. Proofs of mathematical statements are given in Appendices.

\section{Derivative-Based Regularization}
\label{sec:HOR}
The most commonly-used regularizers can be expressed as
\beqn
\mc{R}\br{f}=\int_{{\W}}\Phi\br{\m{D}{f}\br{\mr}}\mbox{d}\mr\,,
\label{eq:General_regularization_form}
\eeqn
where $f$ is an image, $\W\subset\R^2$, $\m{D}$ is the \emph{regularization operator} (scalar or multi-component) acting on the image, and $\Phi\br{\cdot}$ is a potential function. Typical choices for $\m{D}$ are differential operators such as the Laplacian (scalar operator) or the gradient (vectorial operator), while the potential function $\Phi$ usually involves a norm distance. For many years, the preferred choice for the potential function has been the squared Euclidean norm, because of its mathematical tractability and computational simplicity. However, it is now widely documented that non-quadratic potential functions can lead to improved results; they can be designed to be less sensitive to outliers and therefore provide better edge reconstruction. A typical example is TV, which for smooth images corresponds to the $L_1$ norm of the magnitude of the gradient. 

Our present goal is to introduce new regularization functionals of the form of~\eqref{eq:General_regularization_form} which amounts to specifying some suitable linear operator $\m{D}$, and potential function $\Phi$. To do so, certain requirements should be fulfilled. In particular, following the example of TV, we restrict ourselves to regularization operators that commute with translation and scaling, and potential functions that preserve these properties while introducing additional rotation invariance. 
Our motivation for enforcing these invariances is that, similarly to what is the case in many physical systems, one should opt for reconstruction algorithms that lead to solutions which are not affected by transformations of the coordinate system. An additional desirable requirement is that the regularizers should be convex to ensure that if a minimum exists, then this is a global one. Furthermore, convexity permits the design of efficient minimization techniques.

\subsection{Gradient Norm Regularization}
\label{sec:First_Order_Reg}
We would like our regularization operator to be translation and scale invariant. Therefore, a reasonable choice for  $\m{D}$ is some form of derivative operator. Based on this, we first characterize the complete class of gradient-based regularizers satisfying all the required invariances. This is accomplished by Theorem~\ref{thm:grad_norm} which specifies the valid form for the potential functions $\Phi$. The proof of this theorem is given in Appendix~\ref{sec:appendix_A}.
\begin{thm}
Let $\mc{R}\br{f}$ be of the form~\eqref{eq:General_regularization_form}, where $\m{D}$ is the gradient operator and $f$ is continuously differentiable. $\mc{R}\br{f}$ is a translation-, rotation-, and scale-invariant functional, if and only if the potential function $\Phi:\R^2\mapsto\R$ is of the form: $\Phi\br{\grad f\br{\mr}}=c\abs{\grad f\br{\mr}}^\nu$, where $\nu\in\R$ and $c$ is an arbitrary constant. 
\label{thm:grad_norm}
\end{thm}
As a direct consequence of Theorem~\ref{thm:grad_norm}, we see that the following gradient-based regularizers are the only choice of regularization satisfying the required invariance properties; ignoring the multiplicative constant $c$ of the potential function, which can be absorbed by the regularization parameter, we get
\beqn
\mc{R}\br{f}=\int_{{\W}}\abs{\grad f\br{\mr}}^\nu\mbox{d}\mr\,.
\label{eq:Grad_based_regularization_form}
\eeqn
Since we are also interested in convex regularization functionals, we shall focus on cases where $\nu\ge 1$ in~\eqref{eq:Grad_based_regularization_form}. A popular instance of convex functionals arises if we choose $\nu=1$, which corresponds to the TV functional. This regularizer enjoys an additional property, that of contrast covariance.

\subsection{Hessian Schatten-Norm Regularization}
\label{sec:Second_Order_Reg}
As already mentioned, the use of TV, which is the best representative of the gradient-based regularization family, suffers from certain drawbacks. Therefore, for the reasons  specified in the introduction, we are interested in differential operators of higher-order and in particular of the second order. In $N$-dimensions, the complete spectrum of second-order derivatives is embodied in the Hessian operator, 
\beqn
\mc{H}f\br{\mr}=\pr{\begin{array}{cc}
  f_{r_1r_1}\br{\mr} & f_{r_1r_2}\br{\mr} \\
  f_{r_2r_1}\br{\mr} & f_{r_2r_2}\br{\mr} 
  \end{array}},
  \label{eq:Hessian_matrix}
\eeqn
where $f_{r_ir_j}\br{\mr}=\frac{\partial^2}{\partial{r_i}\,\partial{r_j}}f\br{\mr}$. Indeed, with the aid of the Hessian we can compute any second-order derivative of $f\br{\mr}$ as 
$
\mb{D}^2_{\theta,\phi}f\br{\mr}=\m{u}_\t^\transpose\mc{H}f\br{\mr}\m{v}_\phi
$,
where $\m{u}_\t=\br{\cos\t,\,\sin\t}$ and $\m{v}_\phi=\br{\cos\phi,\,\sin\phi}$ are unit-norm vectors specifying the directions of differentiation and $\br{\cdot}^\transpose$ is the transpose operation.

Having specified the regularization operator $\m{D}$, the next step is to investigate which class of potential functions $\Phi$  leads to translation-, rotation-, and scale-invariant second-order regularizers. Next, we provide Theorem~\ref{thm:Hessian_norm} which completely characterizes the form of $\Phi$, under these prerequisites. 
Before presenting this result, we first give the general definition of a Schatten matrix norm~\cite{Bhatia1997}  that will be used in the sequel, and introduce some of the adopted notation. We denote the set of unitary matrices as $\U^n=\cbr{\mX\in\C^{n\times n}:\mX^{-1}=\mX^{\hermitian}}$, where $\C$ is the set of complex numbers and $\br{\cdot}^\hermitian$ is the Hermitian transpose. We also denote the set of positive semidefinite diagonal matrices as $\D^{n_1\times n_2}=\cbr{\mX\in\R_+^{n_1\times n_2}: \mX\br{i,j}=0\,\, \forall\, i\ne j}$, where $\R_+$ is the set of real non-negative numbers.

\begin{definition}[Schatten norms] Let $\mX\in\C^{n_1\times n_2}$ be a matrix with the singular-value decomposition (SVD) $\mX=\mU\mS\mV^\hermitian$, where $\mU \in \U^{n_1}$ and $\mV \in \U^{n_2}$ consist of the singular vectors of $\mX$, and $\mS\in\D^{n_1\times n_2}$ consists of the singular values of $\mX$. The Schatten norm of order $p$ ($\mc{S}_p$ norm) of $\mX$ is defined as
\beqn
\snorm{\mX}{p}=\br{\suml_{k=1}^{\min\br{n_1,n_2}}\ms_k^p\br{\mX}}^{\frac{1}{p}}\,, 
\label{eq:Schatten_def}
\eeqn
where $p \ge 1$, and $\ms_k\br{\mX}$ is the $k$-th  singular value of $\mX$, which corresponds to the $\br{k,k}$ entry of $\mS$.
\label{def:Schatten_norm}
\end{definition}
Definition~\ref{def:Schatten_norm} implies that the $\mc{S}_p$ norm of a matrix $\mX$ corresponds to the $\ell_p$ norm of its singular-values vector $\ms\br{\mX}\in\R_+^{\min\br{n_1,n_2}}$. This further means that all Schatten norms are unitarily invariant. Moreover, we note that the family of $\mc{S}_p$ norms includes three of the most popular matrix norms, i.e., the nuclear/trace norm ($p=1$), the Frobenius norm ($p=2$) and the spectral/operator norm ($p=\infty$).

\begin{thm}
Let $\mc{R}\br{f}$ be of the form~\eqref{eq:General_regularization_form}, where $\m{D}$ is the Hessian operator and $f$ is twice continuously differentiable. $\mc{R}\br{f}$ is a translation-, rotation-, and scale-invariant functional, if and only if the potential function $\Phi:\R^{2\times 2}\mapsto\R$ is of the form: $\Phi\br{\mc{H} f\br{\mr}}=\Phi_0\br{\ml_f\br{\mr}/\snorm{\mc{H}f\br{\mr}}{p}}\snorm{\mc{H} f\br{\mr}}{p}^{\nu}$, where $\nu\in\R$ and $\Phi_0$ is a zero-degree homogeneous function of the Hessian eigenvalues $\ml_f\br{\mr}$.
\label{thm:Hessian_norm}
\end{thm}
The proof of Theorem~\ref{thm:Hessian_norm} is given in Appendix~\ref{sec:appendix_A}. Now, according to it, the admissible second-order regularizers, with respect to the invariance properties of the coordinate system, are those depending on the Schatten norms of the Hessian. If we set $\Phi_0=1$, we obtain the following regularization family
\beqn
\mc{R}\br{f} = \int_{{\W}}\snorm{\mc{H}f\br{\mr}}{p}^{\nu}\mbox{d}\mathbf{r}\,,\forall p\ge 1\mbox{ and } \nu\in\R\,.
\label{eq:Hessian_Reg_general}
\eeqn
To further ensure convexity we need to impose that $\nu\ge 1$. Finally, for our regularizers to also enjoy the contrast-covariance property (similar to TV), we focus on the case where $\nu=1$. Consequently, we define our proposed family of non-quadratic second-order regularization functionals as
\beqn
\mc{R}\br{f} = \int_{{\W}}\snorm{\mc{H}f\br{\mr}}{p}\mbox{d}\mathbf{r}\,,\forall p\ge 1\,.
\label{eq:Hessian_Reg}
\eeqn

The introduced functionals, depending on the Hessian,  lead to piecewise-linear reconstructions. These reconstructions can better approximate the intensity variations observed in natural images than the piecewise-constant reconstructions provided by TV. Thus, they are able to avoid the staircase effect. Moreover, since the Hessian of $f$ at coordinates $\mr$ is a $2\times2$ symmetric matrix, the SVD in the Schatten norm definition reduces to the spectral decomposition and the singular values correspond to the absolute eigenvalues, which can be computed analytically. Now, if we consider the intensity map of the image as a 3-D  differentiable surface, then the two eigenvalues of the Hessian at coordinates $\mr$ correspond to the principal curvatures. They can be used to measure how this surface bends by different amounts in different directions at that point. Therefore, the proposed potential functions, which depend upon those, can be interpreted as scalar measurements of the curvature at a local surface patch.  For example, the $\mc{S}_2$ norm (Frobenius norm) of the Hessian is a scalar curvature index, commonly used in differential geometry, which quantifies lack of flatness of the surface at a specific point. Therefore, we can safely state that the proposed regularizers incorporate curvature information about the image intensity.

Finally, we note that the regularizers obtained for two choices of  $p=2,\infty$, coincide with functionals we considered in our previous work in~\cite{Lefkimmiatis2012}, where we followed another path for extending TV based on rotational averages of directional derivatives. To the best of our knowledge, the Hessian Schatten norms for $p\ne 2,\infty$ have not been considered before in the context of inverse problems.

\section{Variational Image Reconstruction}
\label{sec:image_reconstruction}
From now on, we focus on the discrete formulation of image reconstruction. Hereafter,  to avoid any confusion between the continuous and the discrete domains we will use bold-faced symbols to refer to the discrete manipulation of the problem.
\vspace{-.05cm}
\subsection{Discrete Problem Formulation}
Our approach for reconstructing the underlying image from the measurements is based on the linear observation model
\beqn
\my=\mA\mx+\m{w}\,,
\label{eq:observation_model}
\eeqn
where $\mA \in \R^{M\times N}$ is a matrix that models the spatial response of the imaging device, while $\my \in \R^{M}$ and $\mx\in \R^{N}$ are the vectorized versions of the observed image and the image to be estimated, respectively. Apart from the effect of the operator $\mA$ acting on the underlying image, another perturbation is the measurement noise, which is intrinsic in the detection process. This degradation factor is represented in our observation model by $\m{w}$ that we, here on, will assume to be i.i.d Gaussian noise with variance $\sigma_w^2$.

The recovery of $\mx$ from the measurements $\my$ belongs to the category of linear inverse problems. Usually, for the cases of practical interest, it is  \textit{ill-posed}~\cite{Hansen2006}: the operator $\mA$ is either ill-conditioned or singular. This is dealt with in the variational framework by forming an objective function
\beqn
\varphi\br{\mx}=\frac{1}{2}\norm{\my-\mA\mx}{2}^2+\tau{\psi}\br{\mx}\,,
\label{eq:objective_fun}
\eeqn
whose role is to quantify the quality of a given estimate. The first term, also known as \emph{data fidelity}, corresponds to the negative Gaussian log-likelihood and measures how well a candidate estimate explains the observed data. The second term (\emph{regularization}) encodes our beliefs about certain characteristics of the underlying image. Its role is to narrow down the set of plausible solutions by penalizing those that do not satisfy the assumed properties. The parameter $\tau\ge 0$ provides a balance between the contribution of the two terms. The image reconstruction problem is then cast as the minimization of~\eqref{eq:objective_fun} and leads to a penalized least-squares solution.
\vspace{-.05cm}
\subsection{Discrete Hessian Operator and Basic Notations}
In this work we focus on the class of Hessian Schatten-norm regularizers presented in~\eqref{eq:Hessian_Reg}. In the sequel, we use  
$\dH$ to refer to the discrete version of the Hessian operator. To simplify our analysis, we assume that the image intensities on a $N_x\times N_y$ grid are rasterized in a vector $\mx$ of size $N=N_x\cdot N_y$, so that the pixel at coordinates $\br{i\,,j}$ maps to the $n$th entry of $\mx$ with $n=jN_x+(i+1)$. In this case, the discrete Hessian operator is a mapping $\dH:\R^N\mapsto \mc{X}$, where $\mc{X}=\R^{N\times 2 \times 2}$. For $\mx\in\R^N$, $\dH{\mx}$ is given as
\beqn
\pr{\dH{\mx}}_{n}=\l[\begin{array}{cc}
  \pr{\Dxx\mx}_{n} &   \pr{\Dxy\mx}_{n} \\
  \pr{\Dxy\mx}_{n} &  \pr{\Dyy\mx}_{n}
\end{array}\r]\,,
\eeqn
where $\pr{\cdot}_n$ denotes the $n$th element of the argument, $n=1\,,\ldots,N$, and $\Dxx$, $\Dyy$, and $\Dxy$  denote the forward finite-difference operators~\cite{Denis1983} that approximate the second-order partial derivatives along the two dimensions of the image. If we assume Neumann boundary conditions and use the standard representation of the image rather than the vectorized one, these operators are defined as
\begin{subequations}
\begin{align}
\pr{\Dxx\mx}_{i,j}&=\begin{cases}
 x_{i+2,j}-2 x_{i+1,j}+x_{i,j}, & 1\le i\le N_x-2,\\
 x_{N_x-1,j}-x_{N_x,j},& i\ge N_x-1,\\
\end{cases} \\
\pr{\Dyy\mx}_{i,j}&=\begin{cases}
 x_{i,j+2}-2 x_{i,j+1}+x_{i,j}, & 1\le j\le N_y-2,\\
 x_{i,N_y-1}-x_{i,N_y}, & j\ge N_y-1,\\
\end{cases} \\
\pr{\Dxy\mx}_{i,j}&=\begin{cases}
 x_{i+1,j+1}-x_{i+1,j}-x_{i,j+1}+x_{i,j}, & \\ \quad\quad  1\le i\le N_x-1 \,\,\mb{and}\,\,1\le j\le N_y-1,\\
  0, &\hspace{-5.8cm}\mb{otherwise}\,.
  \end{cases}
  \end{align}
\label{eq:differential_operators}
\end{subequations}

We equip the space $\mc{X}$ with the inner product $\ip{\cdot}{\cdot}_{\mc{X}}$ and norm $\norm{\cdot}{\mc{X}}$. To define them, let $\mX,\mY\in\mc{X}$, with $\mX_n,\mY_n\in\R^{2\times 2}\, \forall\,\, n=1,\ldots,N$. Then we have
\beqn
\ip{\mX}{\mY}_{\mc{X}}=\sum_{n=1}^N\tr{\mY_n^\transpose\mX_n}
\label{eq:inner_product_of_X}
\eeqn
and
\beqn
\norm{\mX}{\mc{X}}=\sqrt{\ip{\mX}{\mX}_{\mc{X}}}\,,
\label{eq:norm_of_X}
\eeqn
where $\tr{\cdot}$ is the trace operator. For the Euclidean space $\R^N$ we use the standard definition of the inner product and of the norm. We denote them by $\ip{\cdot}{\cdot}_{2}$ and $\norm{\cdot}{2}$, respectively.

The adjoint of $\dH$ is the discrete operator $\dH^*:\mc{X}\mapsto \R^N$ such that 
\beqn
\ip{\mY}{\dH\mx}_{\mc{X}}=\ip{\dH^*\mY}{\mx}_{2}\,.
\label{eq:inner_product_relation}
\eeqn
This definition of the adjoint operator is a generalization of the Hermitian transpose for matrices. 
Based on the relation of the inner products in~\eqref{eq:inner_product_relation}, we show in Appendix~\ref{sec:appendix_B} that for any $\mY\in\mc{X}$, it holds that
  \begin{align}
\pr{\dH^*\mY}_n=&\pr{\Dxx^*\mY^{\br{1,1}}}_n+\pr{\Dyy^*\mY^{\br{2,2}}}_n\nonumber\\
&+\pr{\Dxy^*\br{\mY^{\br{1,2}}
+\mY^{\br{2,1}}}}_n\,,
\label{eq:Hessian_adjoin_Op}
\end{align}
 where $\mY^{\br{i,j}}_n$ is the $\br{i,j}$ entry of the $2\times 2$ matrix $\mY_n$ and $\Dxx^*$, $\Dxy^*$, $\Dyy^*$ are the adjoint operators that correspond to backward difference operators with Neumann boundary conditions.

\subsection{Majorization-Minimization Algorithm}
\label{sec:MM}
Next, we present a general method to compute the minimizer of the functional in~\eqref{eq:objective_fun}, under any Hessian-based $\mc{S}_p$ norm regularizer. 
Since these regularizers are non-smooth, our algorithm is based on a majorization-minimization (MM) approach~(cf.~\cite{Daubechies2004}--\cite{Beck2009} for instance).
Under this framework, instead of directly minimizing \eqref{eq:objective_fun}, we find the solution via the successive minimization of a sequence of surrogate functions that upper bound the initial objective function~\cite{Hunter2004}. Our motivation for taking this path is that each of the surrogate functions is simpler to minimize, and we can rely on a gradient scheme that efficiently copes with the large dimensionality of the problem.

To obtain the surrogate functions, we  upper bound the data term of our objective function using the following majorizer~\cite{Daubechies2004,Zibulevsky2010}
\beqn
g\br{\mx,\mx_0}=\frac{1}{2}\norm{\my-\mA\mx}{2}^2+d\br{\mx,\mx_0}\,,
\label{eq:data_majorizer}
\eeqn
where
$
d\br{\mx,\mx_0}=\frac{1}{2}\br{\mx-\mx_0}^\transpose\l[\a \m{I}-\mA^\transpose\mA\r]\br{\mx-\mx_0}
$
is a function that measures the distance between $\mx$ and $\mx_0$.
To come up with a valid majorizer we need to ensure that $d\br{\mx,\mx_0} \ge 0,\,\forall \mx$, with equality if and only if $\mx=\mx_0$. This prerequisite is true if $\a \m{I}-\mA^\transpose\mA$ is positive definite, which implies that $\a>\norm{\mA^\transpose\mA}{}$. The upper-bounded version of the overall objective \eqref{eq:objective_fun} can be written as
\begin{align}
\tilde{\varphi}\br{\mx,\mx_0}
&=\frac{\a}{2}\norm{\mx-\m{z}}{2}^2+\tau\psi\br{\mx}+c\,,
\label{eq:upper_bound}
\end{align}
where $c$ is a constant and $\m{z}=\mx_0+\a^{-1}\mA^\transpose\br{\my-\mA\mx_0}$.
Then, the next step is to iteratively minimize \eqref{eq:upper_bound} w.r.t $\mx$, setting
$\mx_0$ to the previous iteration's solution. As we see, in~\eqref{eq:upper_bound} there is no coupling between $\mx$ and the operator $\mA$ anymore, which turns the minimization task into a much simpler one. In fact, the minimizer of
\eqref{eq:upper_bound} can also be interpreted as the solution of a denoising problem with $\mz$ being the noisy measurements.

\subsection{Proximal Map Evaluation and Matrix Projections}
\label{sec:prox_map}
The MM formulation of our problem relies on the solution of a simpler problem of the form
\beqn
\hat{\mx}=\argmin_{\mx \in \R^{n}}\frac{1}{2}\norm{\mx-\mz}{2}^2+\tau\psi\br{\mx}+\iota_\mc{C}\br{\mx}\,,
\label{eq:prox_min}
\eeqn
where $\iota_\mc{C}$ is the indicator function of  a convex set $\mc{C}$ that represents additional constraints on the solution, such as positivity or box constraints. The convention is that $\iota_\mc{C}\br{\mx}$ takes the value $0$ for $\mx\in\mc{C}$ and $\infty$ otherwise.  If $\vartheta\br{\mx}=\tau\psi\br{\mx}+\iota_\mc{C}\br{\mx}$ is a proper, closed, convex function, then the solution of \eqref{eq:prox_min} is unique and corresponds to the value of the \emph{Moreau proximity operator}~\cite{Combettes2005}, defined as
\beqn
\mbox{prox}_{\vartheta}\br{\mz}=\argmin_{\mx \in \R^{N}}\frac{1}{2}\norm{\mx-\mz}{2}^2+\vartheta\br{\mx}\,.
\eeqn
The proximal map of $\vartheta\br{\mx}$ cannot always be obtained in closed-form, and this is also the case for the regularizers under study. For this reason,
we next present a primal-dual approach that results in a novel numerical algorithm, which can efficiently compute the solution. 

A fundamental ingredient of our proposed algorithm is the orthogonal projection of matrices onto $\mc{S}_q$ norm balls.  This projection can be performed efficiently based on the following proposition, which provides a direct link between vector projections onto $\ell_q$ norm balls and matrix projections onto $\mc{S}_q$ norm balls. This result is new, to the best of our knowledge, and its proof is provided in Appendix~\ref{sec:appendix_B}.  A relevant result that can be considered as a converse statement of Proposition~\ref{thm:matrix_projection} can be found in~\cite[Theorem A.2]{Lewis2008}.
\begin{proposition}[Schatten Norm Projections]
\label{thm:matrix_projection}
Let $\mY\in\C^{n_1\times n_2}$ with SVD decomposition $\mY=\mU\mS\mV^\hermitian$, where $\mU \in \U^{n_1}$, $\mV \in \U^{n_2}$ and $\mS\in\D^{n_1\times n_2}$. The orthogonal projection of $\mY$ onto the
set $\mc{B}_{\mc{S}_q}=\cbr{\mX\in\C^{n_1\times n_2} : \snorm{\mX}{q}\le \rho}$  is given by
\beqn
\mc{P}_{\mc{B}_{\mc{S}_q}}\br{\mY}=\mU\diag{\mc{P}_{\mc{B}_q}\br{\ms\br{\mY}}}\mV^\hermitian\,,\nonumber
\eeqn
where $\diag{\cdot}$ is the operator that maps a vector to a diagonal matrix and $\mc{P}_{\mc{B}_q}$ is the orthogonal projection onto the $\ell_q$ norm ball $\mc{B}_q=\cbr{\m{v}\in\R^{\min\br{n_1,n_2}}_+ : \norm{\m{v}}{q}\le \rho}$ of radius $\rho$.
\end{proposition}
Based on Proposition~\ref{thm:matrix_projection}, we design an algorithm for the orthogonal projection of a matrix $\mY$ onto the convex set $\mc{B}_{\mc{S}_q}$. Our algorithm consists of three steps: (a) decompose  $\mY$  in its singular vectors and singular values by means of the SVD; (b) project its singular values onto the corresponding $\ell_q$ norm ball $\mc{B}_q$; and (c) obtain the projected matrix via singular value reconstruction (SVR) using the projected singular values and the original singular vectors. 

We next describe all of the steps leading to the proposed algorithm that solves the problem
\beqn
\argmin_{\mx \in \mc{C}}\frac{1}{2}\norm{\mx-\mz}{2}^2+\tau\norm{\dH{\mx}}{1,p}\,\, \forall p\ge 1\,.
\label{eq:denoise_opt}
\eeqn
With $\norm{\dH{\mx}}{1,p}$ we denote the discrete version of our proposed regularization family~\eqref{eq:Hessian_Reg}, where \norm{\cdot}{1,p} stands for 
the mixed $\ell_1$-$\mc{S}_p$ norm, which for an argument $\mP=\pr{\mP_1^\transpose, \mP_2^\transpose,\ldots, \mP_N^\transpose }^\transpose \in \mc{X}$ is defined as
\beqn
\norm{\mP}{1,p}=\suml_{n=1}^{N}\snorm{\mP_{n}}{p}\,,\forall p\ge 1\,.
\label{eq:mixed_l1_Sp}
\eeqn
The discrete form of our regularizers highlights their relation to the sparsity-promoting group norms, which are commonly met in the context of compressive sensing (see~\cite{Bach2011}, for instance). However, a significant difference is that in our case the mixed norm is a vector-matrix norm rather than a vector-vector norm. Therefore, while the machinery we are using shares some similarities with the one employed in the group vector-norm case, there are important differences, with the most pronounced being the projection step.

Since the operator of our choice is the Hessian, which produces $2\times 2$ symmetric matrices at every coordinate of $\mx$, $\mP_{n}\in\mathbb{S}^{2}$ in~\eqref{eq:mixed_l1_Sp}, where $\mathbb{S}^2=\cbr{\mX\in \R^{2\times 2}: \mX^\transpose=\mX}$.  However, for reasons of completeness, in the following lemma, where we derive the dual of the $\ell_1$-$\mc{S}_p$ norm, we consider the more general case $\mP_n\in\C^{n_1\times n_2}$. The proof of Lemma~\ref{lemma:lemma_dual_norm} is provided in Appendix~\ref{sec:appendix_B} and
follows a similar line of thought with the one presented in~\cite[Lemma 1]{Sra2011}. The latter is about the dual norm of a mixed $\ell_1$-$\ell_p$ vector norm. 

\begin{lemma}
\label{lemma:lemma_dual_norm}
Let $p\ge 1$, and let $q$ be the conjugate exponent of $p$,  i.e., $\frac{1}{p}+\frac{1}{q}=1$. Then, the mixed norm $\norm{\cdot}{\infty,q}$ is dual to the mixed norm $\norm{\cdot}{1,p}$.
\end{lemma}
Using Lemma \ref{lemma:lemma_dual_norm} and noting that the dual of the dual norm is the original norm~\cite{Rockafellar1970}, we  write~\eqref{eq:mixed_l1_Sp} in the equivalent form
\beqn
\norm{\mP}{1,p}=\max_{\mW\in\mc{B}_{\infty,q}}\ip{\mW}{\mP}_{\mc{X}}\,,
\label{eq:dual_mixed_l1_Sp}
\eeqn
where $\mc{B}_{\infty,q}$ denotes the $\ell_\infty$-$\mc{S}_q$ unit-norm ball, defined as 
\begin{align}
\mc{B}_{\infty,q}=\Big\{\mW&=\pr{\mW_1^\transpose, \mW_2^\transpose,\ldots, \mW_N^\transpose }^\transpose \in \mc{X}: \nonumber\\ &\snorm{\mW_n}{q}\le 1, \forall n=1,\ldots,N\Big\}.
\label{eq:unit_norm_ball}
\end{align}
This alternative definition of the mixed $\ell_1$-$\mc{S}_p$ norm allow us to express it in terms of an inner product that involves the dual variable $\mW$ and the unit-norm ball $\mc{B}_{\infty,q}$. Moreover, from~\eqref{eq:unit_norm_ball} it is straightforward to see that the orthogonal projection onto $\mc{B}_{\infty,q}$ is obtained by projecting separately each submatrix $\mW_n$ onto a unit-norm $\mc{S}_q$ ball ($\mc{B}_{\mc{S}_q}$).

Using~\eqref{eq:dual_mixed_l1_Sp} we re-write~\eqref{eq:denoise_opt} as
\beqn
\hat{\mx}=\argmin_{\mx \in \mc{C}}\frac{1}{2}\norm{\mx-\mz}{2}^2+\tau\max_{\mW\in\mc{B}_{\infty,q}}\ip{\mW}{\dH{\mx}}_\mc{X}\,.
\eeqn
This formulation naturally leads us to the following minimax problem
\beqn
\min_{\mx \in \mc{C}}\max_{\mW\in\mc{B}_{\infty,q}}\mc{L}\br{\mx,\mW}\,,
\label{eq:minimax}
\eeqn
where
\begin{align}
\mc{L}\br{\mx,\mW}
&=\frac{1}{2}\norm{\mx-\mz}{2}^2+\tau\ip{\dH^*\mW}{\mx}_2\,.
\end{align}
Since the function $\mc{L}\br{\mx,\mW}$ is strictly convex in $\mx$ and concave in $\mW$, we have the guarantee that a saddle-value is attained \cite{Rockafellar1970}, and, thus, the order of the minimum and the maximum in \eqref{eq:minimax} does not affect the solution. This means that there exists a saddle-point $\br{\hat\mx,\hat\mW}$ that leads to a common value when the minimum and the maximum are interchanged, i.e.,
\beqn
\min_{\mx \in \mc{C}}\max_{\mW\in\mc{B}_{\infty,q}}\mc{L}\br{\mx,\mW}=
\mc{L}\br{\hat{\mx},\hat{\mW}}=
\max_{\mW\in\mc{B}_{\infty,q}}\min_{\mx \in\mc{C}}\mc{L}\br{\mx,\mW}\,.
\label{eq:minmax_maxmin_equiv}
\eeqn
Based on this observation, we can now define the primal and dual problems by identifying the primal and dual objective functions, respectively. The l.h.s of \eqref{eq:minmax_maxmin_equiv} corresponds to the minimization of the primal objective function $\varrho\br{\mx}$, and the r.h.s to the maximization of the dual objective function $s\br{\mW}$,
\beqn
\varrho\br{\mx}=\max_{\mW\in\mc{B}_{\infty,q}}\mc{L}\br{\mx,\mW}=\frac{1}{2}\norm{\mx-\mz}{2}^2+\tau\norm{\dH{\mx}}{1,p}\,,
\eeqn
\begin{align}
s\br{\mW}=&\min_{\mx\in\mc{C}}\mc{L}\br{\mx,\mW}\nonumber\\
=&\frac{1}{2}\norm{\mc{P}_\mc{C}\br{\m{v}}-\m{v}}{2}^2
+\frac{1}{2}\norm{\mz}{2}^2-\frac{1}{2}\norm{\m{v}}{2}^2\,,
\label{eq:dual_opt}
\end{align}
where $\mc{P}_\mc{C}$ is the orthogonal projection onto the convex set $\mc{C}$ and $\m{v}=\mz-\tau\dH^*\mW$. Therefore, \eqref{eq:minmax_maxmin_equiv} indicates that we can obtain the minimizer $\hat{\mx}$ of $\varrho\br{\mx}$ from the maximizer $\hat{\mW}$ of $s\br{\mW}$ through the relation
\beqn
\hat{\mx}=\mc{P}_\mc{C}\br{\mz-\tau\dH^*\hat{\mW}}\,.
\label{eq:primal_dual_relation}
\eeqn
This last relation is important, since in contrast to the primal problem~\eqref{eq:denoise_opt}, which is not continuously differentiable, the dual one involves the smooth function $s\br{\mW}$. We can therefore solve it by exploiting its gradient. Indeed, using the property that the gradient of a function $h\br{\mx}=\norm{\mx-\mc{P}_{\mc{C}}\br{\mx}}{2}^2$ is well defined and is  equal to $\grad{h}\br{\mx}=2\br{\mx-\mc{P}_{\mc{C}}\br{\mx}}$~\cite[Lemma 4.1]{Beck2009b}, we compute the gradient of $s\br{\mW}$ as
\beqn
\grad{s}\br{\mW}=\tau\dH\mc{P}_{\mc{C}}\br{\mz-\tau\dH^*\mW}.
\label{eq:grad_of_s}
\eeqn
Therefore, the solution of our primal problem~\eqref{eq:denoise_opt} is obtained in two steps: (a) we find the maximizer of the dual objective function~\eqref{eq:dual_opt} as described next, and (b) we obtain the solution through~\eqref{eq:primal_dual_relation}.

\subsection{Maximization of the Dual Objective}
\label{sec:max_dual}
At this point, a main issue we need to deal with, is that the Hessian operator $\dH$ does not have an empty null space and, thus, a stable inverse does not exist. Consequently, we cannot opt for a closed-form solution for the maximizer $\hat{\mW}$ of $s\br{\mW}$. This means that we have to resort to a numerical iterative scheme. In this work, we employ Nesterov's iterative method~\cite{Nesterov1983} for smooth functions. This is a gradient-based scheme that exhibits convergence rates of one order higher than the standard gradient-ascent method. To ensure convergence of the algorithm, we need to choose an appropriate step-size. Since our dual objective is smooth with Lipschitz continuous gradient,  we can use a constant step-size, thus, avoid a line search at every iteration. An appropriate step-size is equal to the inverse of the Lipschitz constant of $\grad{s}\br{\mW}$.  We derive an upper bound of this Lipschitz constant in the following proposition, whose proof is given in Appendix~\ref{sec:appendix_B}.

\begin{proposition}
\label{prop:Lipschitz_constant}
Let $L\br{s}$ denote the Lipschitz constant of $\grad{s}\br{\mW}$ of the dual objective function defined in~\eqref{eq:dual_opt}. Then, it holds that
\beqn
L\br{s}\le 64\tau^2\,.
\eeqn
\end{proposition}

From~\eqref{eq:minmax_maxmin_equiv} and~\eqref{eq:dual_opt} it is clear that the maximizer of our dual objective can be derived by solving the constrained maximization problem
\begin{align}
\hat{\mW}=&\argmax_{\mW\in\mc{B}_{\infty,q}}\frac{1}{2}\norm{\mc{P}_\mc{C}\br{\m{v}}-\m{v}}{2}^2
-\frac{1}{2}\norm{\m{v}}{2}^2,
\end{align}
with $\m{v}=\mz-\tau\dH^*\mW$.
A necessary step towards this direction is to compute the projection onto the set $\mc{B}_{\infty,q}$, defined in \eqref{eq:unit_norm_ball}. This operation is accomplished by projecting independently each of the $N$ components $\mW_n$ of $\mW$ onto the set $\mc{B}_{\mc{S}_q}=\cbr{\mX\in\S^2 : \snorm{\mX}{q}\le 1}$. This projection is performed efficiently following the three steps of the algorithm we proposed in Section~\ref{sec:prox_map}, which is based on our Proposition~\ref{thm:matrix_projection}.

Steps (a) and (c) are fairly easy to implement. Specifically, since the matrices of interest are $2\times 2$ symmetric, we compute the SVD and the SVR  steps in closed-form. Then, the most cumbersome part of our algorithm is the $\ell_q$ norm projection of the singular values, which for general values of $q$ does not exist in closed form. Fortunately, this operation is still feasible thanks to the recently developed $\ell_q$ norm projection algorithm~\cite{Sra2011}. This projection method is based on an efficient proximity algorithm for $\ell_q$ norms~\cite{Liu2010}. Moreover, in Section~\ref{sec:Closed_Form_Projections}  we report three cases of $\mc{S}_q$ norms, $q=1,2,\infty$, where the projection can be evaluated in closed-form.

\subsection{Closed Form of $\mc{S}_q$-Norm Projections for $q=1,2,\infty$.}
\label{sec:Closed_Form_Projections}
From Proposition~\ref{thm:matrix_projection}, we know that the matrix projection onto the $\mc{B}_{\mc{S}_2}$ unit-norm ball is associated with the projection of the singular values of the matrix onto the $\mc{B}_2$ ball. The latter is computed by normalizing the elements of the corresponding vector by their Euclidean norm. Therefore, we have that
\beqn
\mc{P}_{\mc{B}_{\mc{S}_2}}\br{\mW_n}=\begin{cases}
\frac{\mW_n}{\norm{\mW_n}{F}} &, \mbox{if $\norm{\mW_n}{F}>1$}\\
\mW_n &, \mbox{if $\norm{\mW_n}{F}\le1$}\,.
\end{cases}
\eeqn
This situation is advantageous since it allows us to avoid both the SVD and the SVR steps. Consequently, this drastically reduces the complexity of computing the projection.
To compute the projection onto $\mc{B}_{\mc{S}_\infty}$, we use that the projection onto the  $\mc{B}_{\infty}$ unit-norm ball corresponds to setting the elements that have an absolute value greater than one to one, and adding back their original sign. Therefore, we readily get 
\beqn
\mc{P}_{\mc{B}_{\mc{S}_\infty}}\br{\mW_n}=\mU\diag{\min\br{\ms\br{\mW_n},\m{1}}}\mV^\hermitian\,,
\eeqn
where $\m{1}$ is a vector with all elements set to one and the $\min$ operator is applied component-wise.
Note that, this result is directly related to the singular value thresholding (SVT) method~\cite{Cai2008,Ma2011} developed in the field of matrix rank minimization. The derivation of SVT in~\cite{Cai2008,Ma2011} is technical. It relies on the characterization of the subgradient of the nuclear norm~\cite{Watson1992}. By contrast, in our case the result comes out naturally as an immediate consequence of Proposition~\ref{thm:matrix_projection} and of the duality between the spectral and nuclear matrix norms.
Finally, the projection of a matrix onto the $\mc{B}_{\mc{S}_{1}}$ unit-norm ball is related to the projection of its singular values onto the $\mc{B}_1$ unit-norm ball. The latter is computed by the soft-thresholding operator $S_\gamma\br{\ms\br{\mW_n}}=\max\br{\ms\br{\mW_n}-\gamma,0}$~\cite{Donoho1995}, where the  $\max$ operator is applied component-wise. Therefore, based on Proposition~\ref{thm:matrix_projection}, we have that
\beqn
\mc{P}_{\mc{B}_{\mc{S}_1}}\br{\mW_n}=\mU\diag{S_\gamma\br{\ms\br{\mW_n}}}\mV^\hermitian\,.
\eeqn
This last projection cannot in general be computed in closed form. The reason is that the threshold $\gamma$ is not known in advance and needs to be estimated. This can be accomplised using one of the existing methods available in the literature~\cite{Candes2005}--\cite{Duchi2008}.
Fortunately, in our case, the singular vectors are of low dimensionality, $\ms\br{\mW_n}\in\R_+^{2}$ and, thus, we derive $\gamma$ analytically, as
\begin{align}
\gamma=\begin{cases}
0 &, \mb{if $\ms_1\br{\mW_n}\le1-\ms_2\br{\mW_n}$}\,,\\
\frac{\ms_1\br{\mW_n}+\ms_2\br{\mW_n}-1}{2} &, \mb{if $1-\ms_2\br{\mW_n}<\ms_1\br{\mW_n}\le1+\ms_2\br{\mW_n}$}\,,\\
\ms_1\br{\mW_n}-1 &, \mb{if $\ms_1\br{\mW_n}>1+\ms_2\br{\mW_n}$}\,,
\end{cases}
\label{eq:threshold_rule}
\end{align}
where  the singular values are sorted in a decreasing order, i.e., $\ms_1\br{\mW_n}\ge\ms_2\br{\mW_n}$.

\subsection{Numerical Algorithm}
Equipped with all the necessary ingredients, we conclude with a summarized description of the complete optimization algorithm. Our method consists of two components that interact. The first component is responsible for the majorization of the objective function, as we described in  Sec.~\ref{sec:MM}, while the second one undertakes the  minimization of the resulting upper-bounded version. Then, the algorithm proceeds by iteratively minimizing the majorizer that is formed based on the solution of the previous iteration. Since the convergence of this scheme can be slow in practice, to speed it up we employ the FISTA algorithm~\cite{Beck2009}. This method exhibits state-of-the-art convergence rates by combining  two consecutive iterates, in an optimum way.  A description of our image reconstruction approach that is based on the monotone version of FISTA (MFISTA)~\cite{Beck2009b} is given in Algorithm 1.  The sub-routine \texttt{denoise} corresponds to the second component that finds the solution of~\eqref{eq:denoise_opt}. This minimizer is related to the proximal map $\mbox{prox}_{\tau\norm{\dH\cdot}{1,p}}$ but we can also interpret it as a denoising step under Hessian-based $\ell_1$-$\mc{S}_p$ norm regularization.  The computation of the  $\texttt{denoise}$ sub-routine is described in Algorithm 2 and is based on the primal-dual formulation that we proposed in Secs.~\ref{sec:prox_map} and \ref{sec:max_dual}. 

Finally, regarding the computational complexity of the algorithm, it is only mildly higher than that of TV's. The extra computational cost is due to (a) the use of a tensor (Hessian) instead of a vectorial (gradient) operator and (b) the projections onto the $\mc{B}_{\mc{S}_q}$ balls instead of the $\mc{B}_2$ ball. Our projections are somewhat more expensive because of the SVD and SVR steps. However, these steps are computed in closed form. It is also worth mentioning that the proposed algorithm is highly parallelizable, since all the involved operations are performed independently for each pixel of the image. 
\begin{figure}[t]
\centering
   \includegraphics[scale=.8]{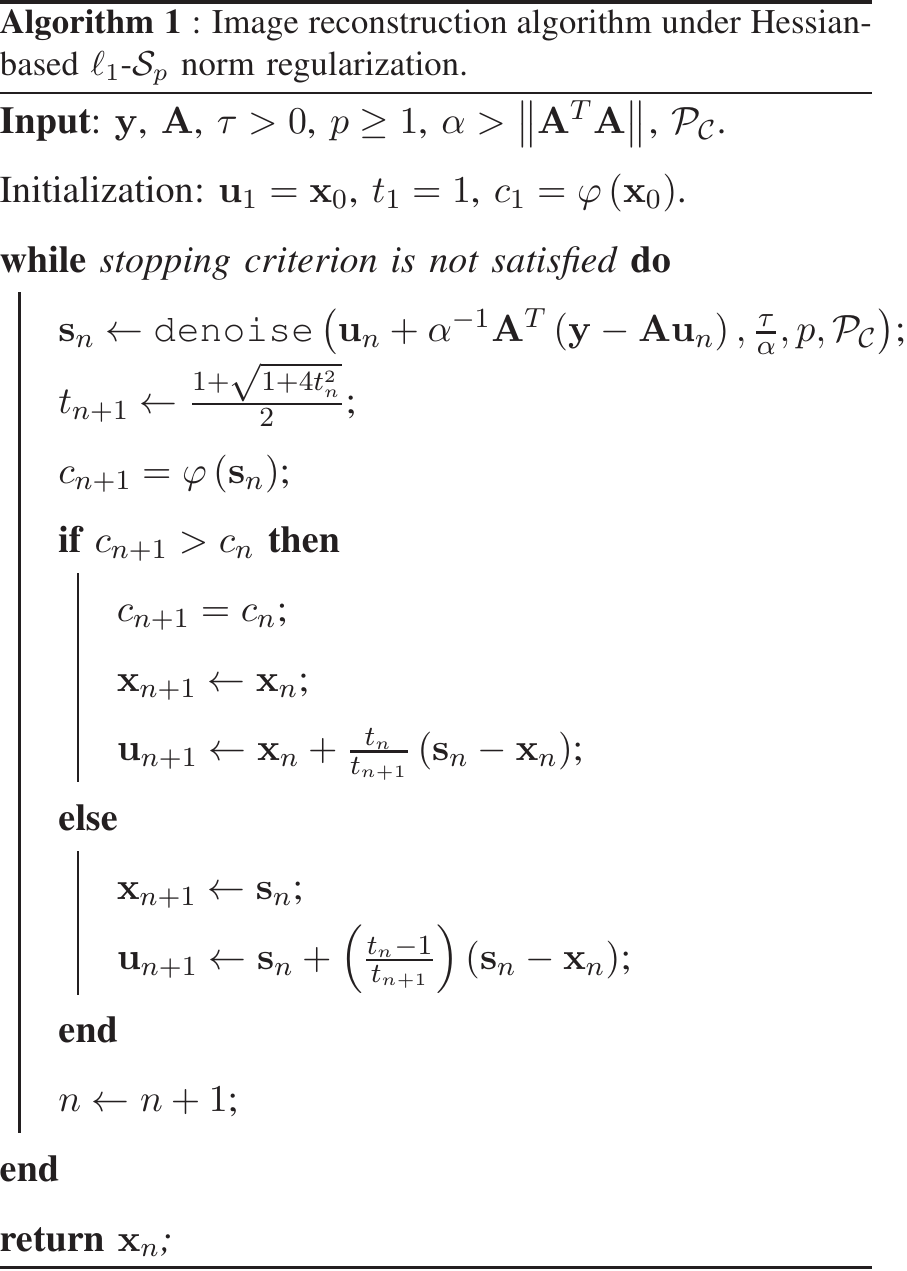}
\label{alg:deconv}
\vspace{-.4cm}
\end{figure}
\begin{figure}[t]
\centering
  \includegraphics[scale=.8]{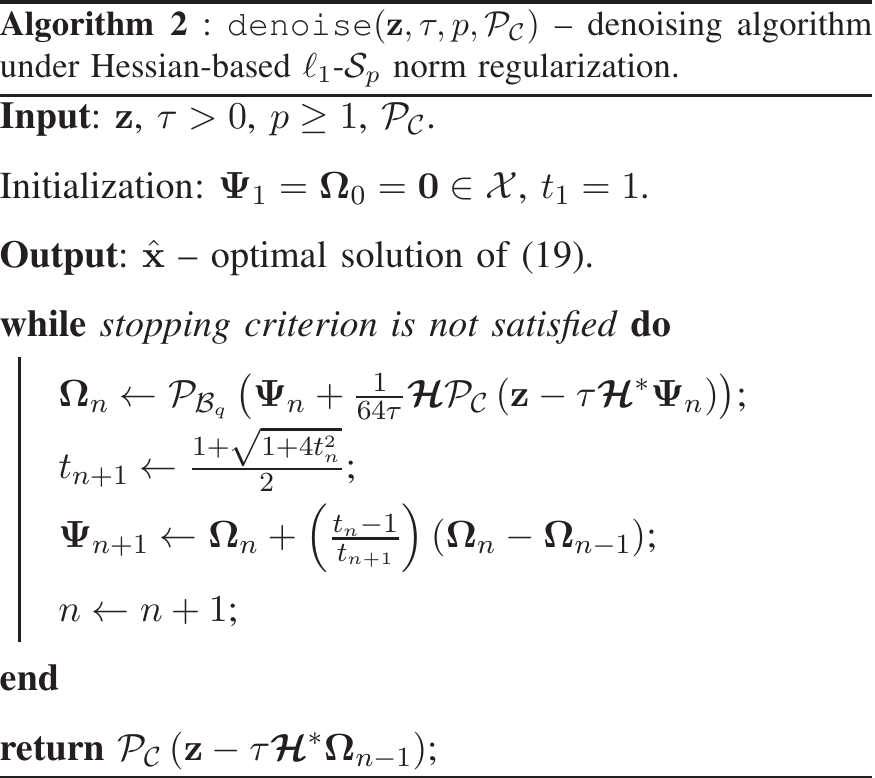}
\label{alg:denoising}
\vspace{-.4cm}
\end{figure}

\section{Experimental Results}
\label{sec:results}
To evaluate the effectiveness of our proposed regularization framework, we report results for several linear inverse imaging problems. In particular, we consider the problems of image deblurring, sparse reconstruction from random samples, image interpolation and image zooming. For the image deblurring problem we compare our results against those obtained by using three alternative methods; namely, TV regularization, regularization with the fully redundant Haar wavelet transform, and the image deblurring version of the BM3D patch-based method~\cite{Dabov2008}.  In Haar's case, we use the frame analysis framework since it has been reported in the literature~(c.f \cite{Selesnick2009}) that the frame synthesis framework usually leads to inferior results. For the rest of the inverse problems we provide comparisons against TV and quadratic derivative-based regularizers. 
\vspace{-.2cm}
\subsection{Restoration Setting}
For the image deblurring experiments, we use a set of 8 grayscale images\footnote{Three of  these images along with the motion-blur kernel used in the experiments were obtained from~\url{http://www.wisdom.weizmann.ac.il/~levina/papers/LevinEtalCVPR09Data.rar}.} which have been normalized so that their intensities lie in the range of $\pr{0\,,1}$. 

The performance of the methods under comparison is assessed for various blurring kernels and different noise levels. In particular, in our experiments we employ three \emph{point spread functions} (PSFs) to produce blurred versions of the images. We use a  Gaussian PSF of standard deviation $\sigma_b=4$, a moving average (uniform) PSF, and a motion-blur kernel. The first two PSFs have a support of $9\times 9$ pixel while the third one has a support of $19\times 19$ pixel. As an additional degradation factor we consider Gaussian noise of three noise levels corresponding to a blurred signal-to-noise-ratio (BSNR) of $\cbr{15, 20, 25}$ dB, respectively. The BSNR is defined as $\mbox{BSNR}=\mbox{var}\br{\mA\mx}/\sigma_w^2$ where $\mbox{var}\br{\mA\mx}$ is the variance of the blurred image and $\sigma_w$ is the standard deviation of the noise.

Regarding the restoration task, for the methods that involve the minimization of an objective function this is performed under the constraint that the restored intensities must lie in the convex set $\mc{C}=\cbr{\mx \in \R^N | x_n\in\pr{0\,,1}\forall n=1,\ldots,N}$. To accomplish that, we use the corresponding projection operation, $\mc{P}_\mc{C}$, which for a vector $\mx$ 
amounts to setting the elements that are less than zero and greater to one, to zero and one, respectively.
For the Hessian-based functionals we use the minimization method proposed in Section~\ref{sec:image_reconstruction}, while for TV- and Haar-based ones we employ the algorithm of~\cite{Beck2009b}. The latter belongs to the same category of minimization algorithms as ours with a comparable convergence behavior. The rationale for this choice is that, the quality of the restoration will not depend on the choice of the minimization strategy but rather on the choice of the regularizer.  In all cases, the stopping criterion is set to either reaching a relative normed difference of $10^{-5}$ between two successive estimates, or a maximum of 100 MFISTA iterations. We also use 10 inner iterations for the solution of the corresponding denoising problem. Moreover, instead of using the true PSF that produces the blurred images, we use a slightly perturbed version by adding Gaussian noise of standard deviation $10^{-3}$. The motivation is to test the performance of the algorithms under more realistic conditions, since, in practice the employed PSF normally contains some error and thus deviates from the true one.
Finally, in all the reported experiments the quality of the reconstruction is evaluated in terms of an increase in the SNR (ISNR), measured in dB. The ISNR is defined as
$\mbox{ISNR}= 10\log_{10}\br{\mbox{MSE}_{\mbox{\scriptsize in}}/{\mbox{MSE}_{\mbox{\scriptsize out}}}}$,
where $\mbox{MSE}_{\mbox{\scriptsize in}}$ and $\mbox{MSE}_{\mbox{\scriptsize out}}$ are the mean-squared errors between the degraded and the original image, and the restored and the original image, respectively.

\vspace{-.3cm}
\subsection{Image Restoration on Standard Test Images}
In Table~\ref{tab:results} we provide comparative restoration results for all the test images and all the combinations of degradation conditions (PSF and noise level). To distinguish between the different Hessian-based regularizers, we refer to them  as $\mc{HS}_k$ with $k$ denoting the order of the Schatten norm. We report the results obtained by using Schatten norms of order one, two and infinity, which correspond to the well-known nuclear, Frobenius and spectral matrix norms, respectively. For the sake of consistency among comparisons, the reported results for each regularizer, including Haar and TV, are derived using the individualized (w.r.t. degradation conditions) regularization parameter $\tau$,  that gives the best ISNR performance. The results of the BM3D algorithm are also optimized by providing the true standard deviation of the Gaussian noise. 

\begin{table*}[!t]
 \centering
  \caption{ISNR comparisons on image restoration for three blurring kernels and three noise levels}
    \includegraphics[scale=.75]{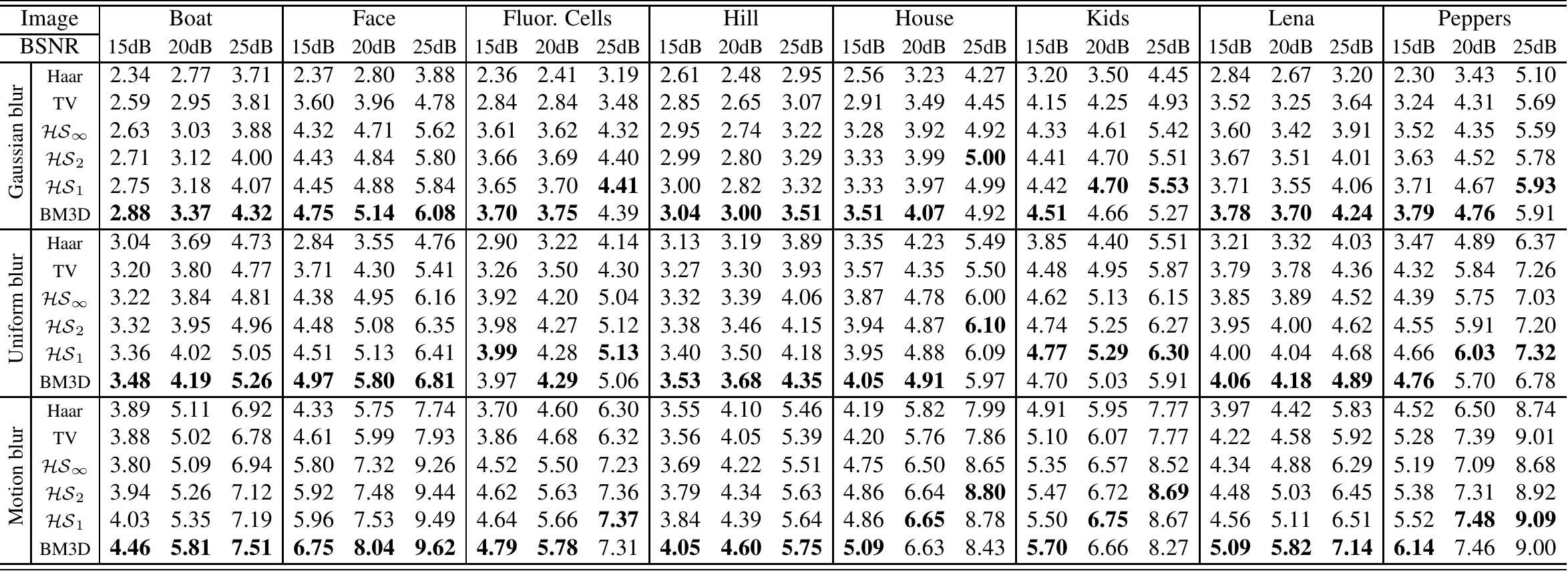}
 \label{tab:results}
\end{table*}

On average, for all the tested images and degradation conditions, the BM3D algorithm produces the best PSNR scores. However, despite its non-adaptive nature, our regularization scheme manages to provide comparable results. Regarding comparisons among the regularization techniques, the Hessian-based framework leads to improved quantitative results compared to those of Haar and TV. The best SNR improvement, on average, is achieved for the $\mc{HS}_1$ regularizer, while comparable results are also obtained for the $\mc{HS}_2$ regularizer. While the $\mc{HS}_\infty$ regularizer  outperforms Haar and TV most of the time, the improvement is less pronounced than that of the other two regularizers. We can thus conclude, that as the order of the Schatten norm moves  closer to 1, the reconstruction results improve. This can be attributed to the fact that, in the extreme case of order infinity, the corresponding regularizer takes into account only the maximum absolute eigenvalue and thus fails to include additional information possibly provided by the second eigenvalue. Overall, the improvement in performance over Haar and TV can be quite substantial (more than 0.5 dB), which justifies Hessian-based regularization as a viable alternative approach.

\begin{figure}[!t]
  \centering
        \includegraphics[scale=.35]{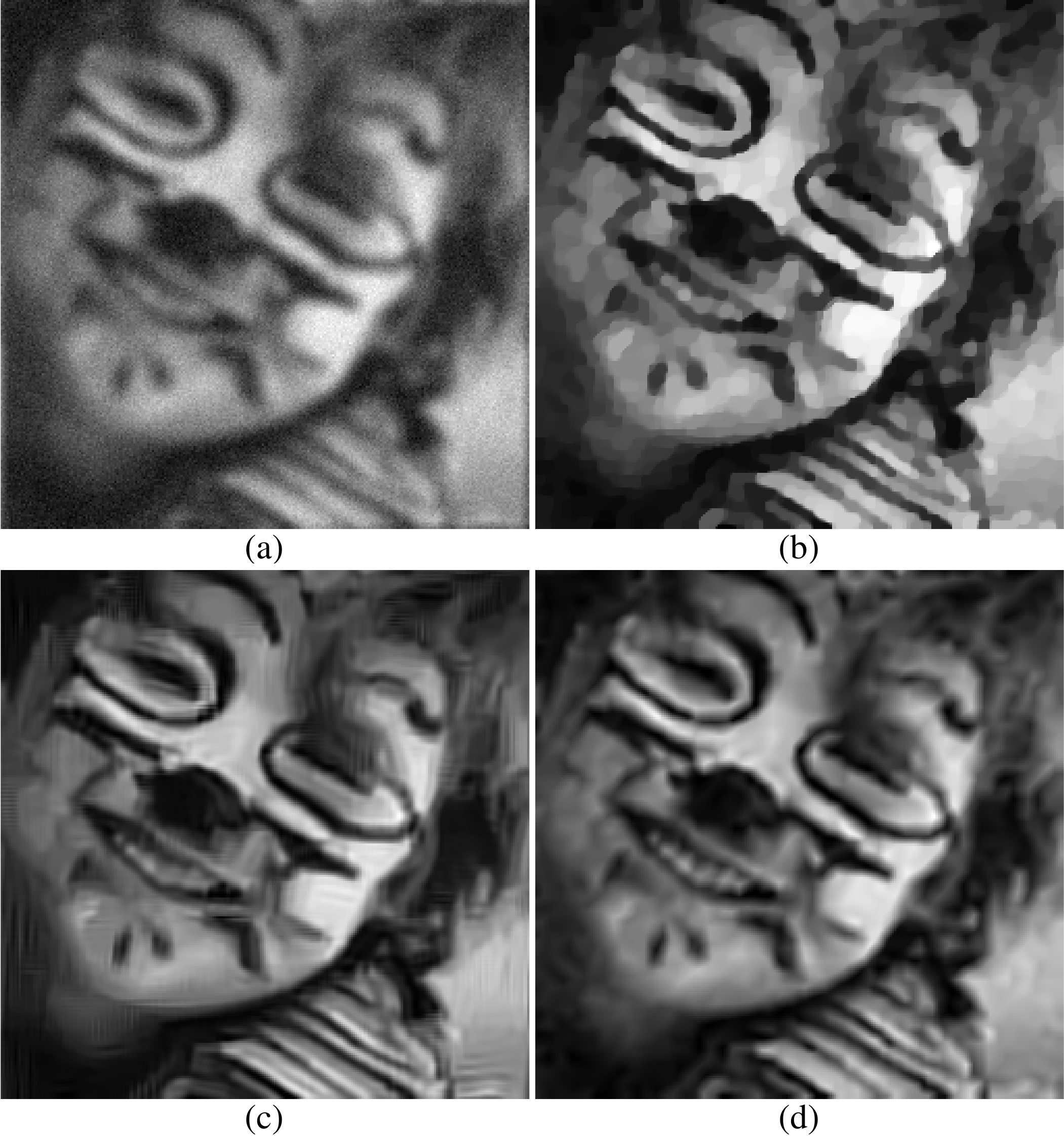}\vspace{-.2cm}
 \caption{Restoration of the Face image degraded by Gaussian blurring and Gaussian noise corresponding to a BSNR level of 15 dB. (a) Degraded image (PSNR = 21.76 dB), (b) TV result (PSNR = 25.36 dB),
 (c) BM3D result (PSNR = 26.51 dB), and (d) $\mc{HS}_1$ result (PSNR = 26.21 dB)}
 \label{fig:face_results}
 \vspace{-0.5cm}
\end{figure}

\begin{figure}[!t]
  \centering
    \includegraphics[scale=.35]{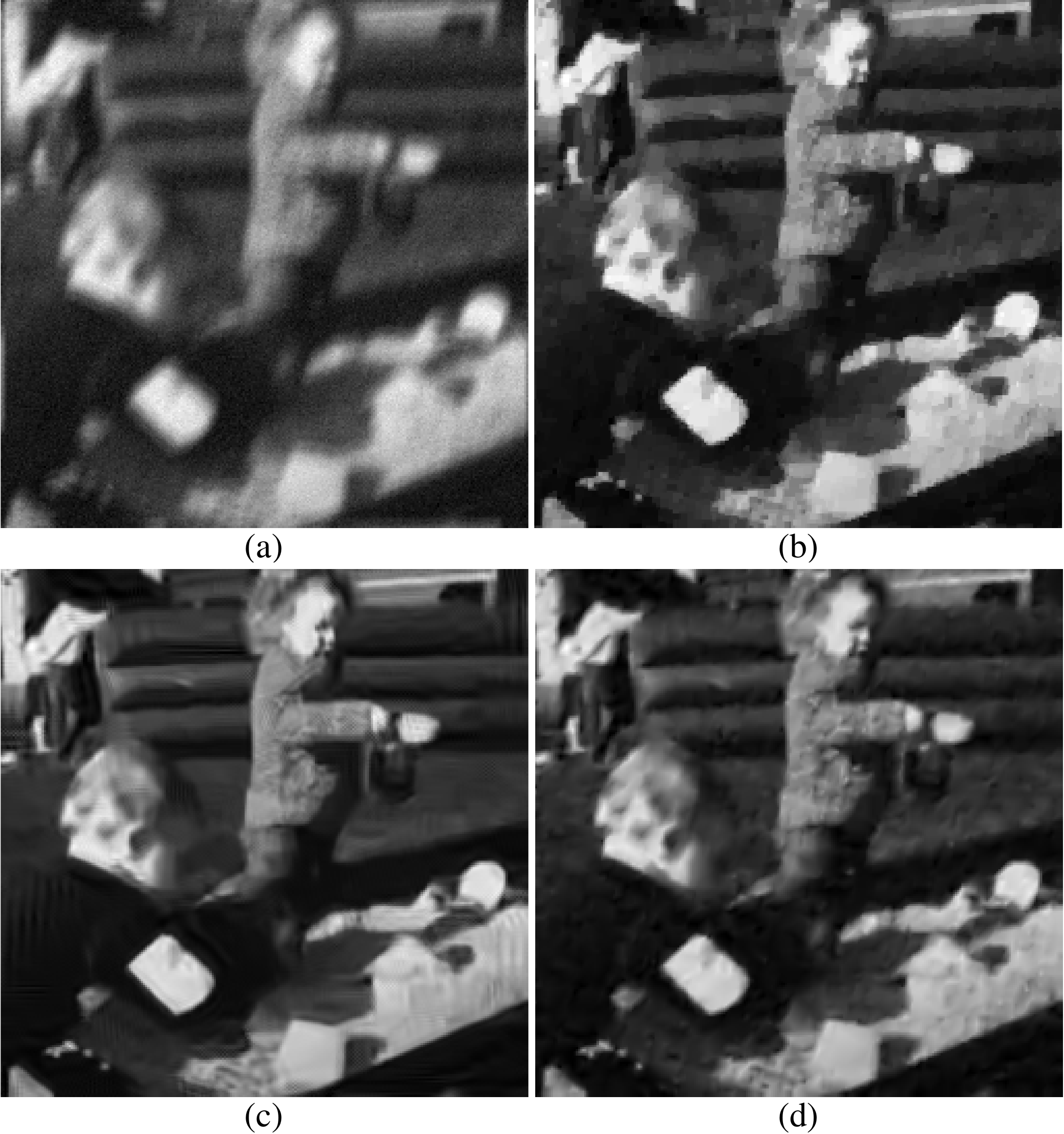}\vspace{-.2cm}
 \caption{Restoration of the Kids image degraded by motion blurring and Gaussian noise corresponding to a BSNR level of 20 dB. (a) Degraded image (PSNR = 21.84 dB), (b) Haar result (PSNR = 27.79 dB),
 (c) BM3D result (PSNR = 28.50 dB), and (d) $\mc{HS}_1$ result (PSNR = 28.59 dB)}
 \label{fig:kids_results}
 \vspace{-.4cm}
\end{figure}

\begin{figure}[!t]
  \centering
    \includegraphics[scale=.35]{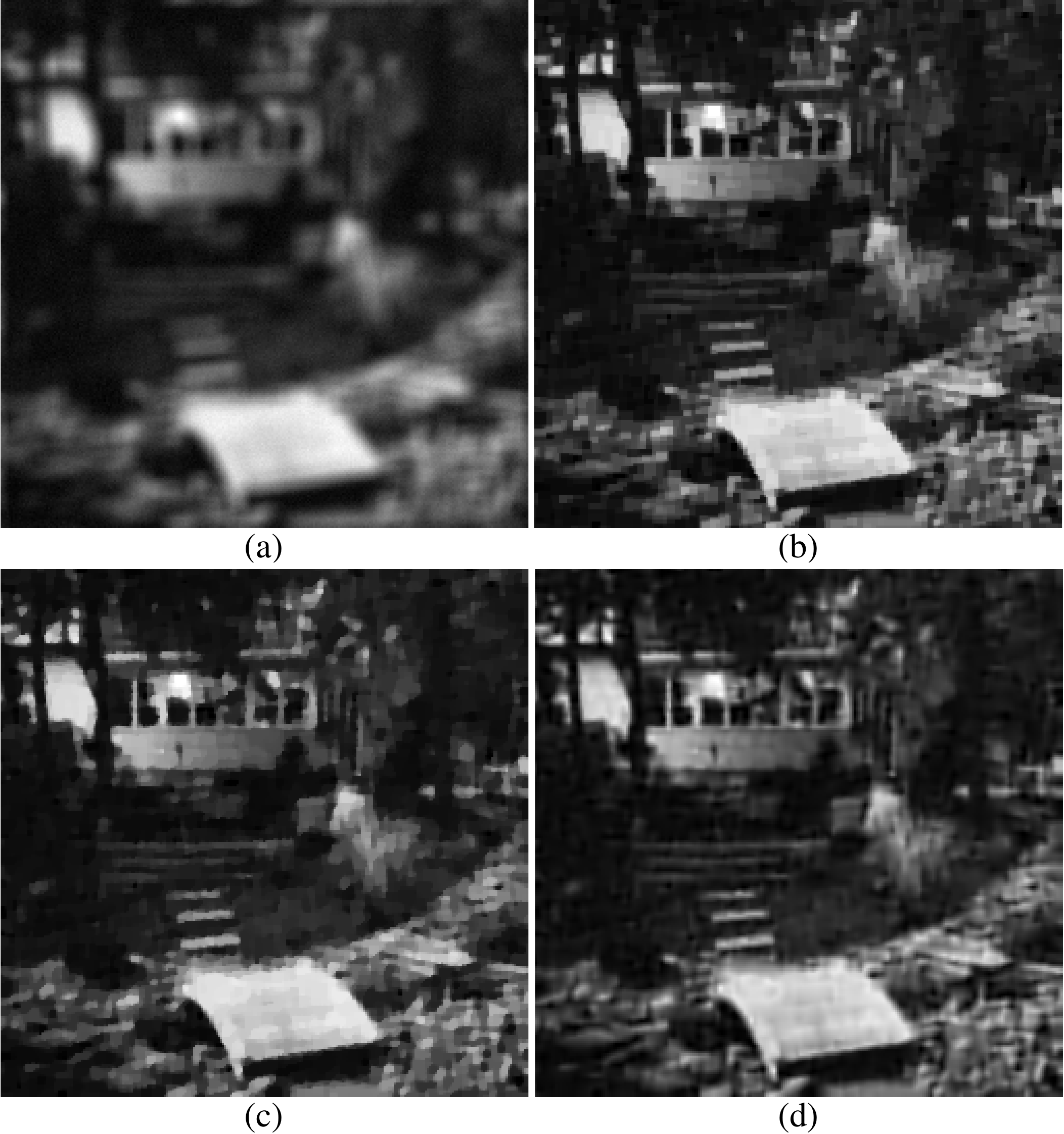}\vspace{-.2cm}
    \caption{Restoration of the House image degraded by uniform blurring and Gaussian noise corresponding to a BSNR level of 25 dB. (a) Degraded image (PSNR = 20.43 dB), (b) Haar result (PSNR = 25.92 dB),
 (c) TV result (PSNR = 25.93 dB), and (d) $\mc{HS}_2$ result (PSNR =  26.53 dB)}
 \label{fig:house_results}
 \vspace{-.5325cm}
\end{figure}

Beyond the ISNR comparisons, the effectiveness of the proposed method can also be visually appreciated by inspecting the representative Face, Kids and House deblurring examples of Figures~\ref{fig:face_results} - \ref{fig:house_results}. From these examples we can verify our initial claims, that TV regularization leads to image reconstructions that suffer from the presence of heavy block artifacts. These artifacts become more evident in regions where the image is characterized by smooth intensity transitions, and they are responsible for shuffling details of the image and broadening its fine structures. See for example the TV solution in Fig.~\ref{fig:face_results}, where the image has cartoon-like appearance.  Similar blocking effects, which are even more pronounced, appear on the Haar-based reconstructions. On the other hand, even in cases where the presence of the noise is significant,  the Hessian-based regularizers manage to avoid introducing pronounced artifacts and thus, they lead to reconstructions that are more faithful representations of the original content of the image. Comparing our results with those of BM3D, we note that even in cases where the final PSNR favors the latter reconstruction, such as in Fig.~\ref{fig:face_results}, our restored images have certain advantages. For example, by a careful inspection of Figs.~\ref{fig:face_results} and \ref{fig:kids_results}, one can clearly observe the presence of ripple-like artifacts in the BM3D solutions which do not appear in the Hessian-based reconstructions. 


\subsection{Deblurring of Biomedical Images}
Our interest in image restoration is mostly motivated by the problem of microscopy image deblurring.  In widefield microscopy, the acquired images are degraded by out-of-focus blur due to the poor localization of the microscope's PSF. This severely reduces our ability to clearly distinguish fine specimen structures. Since a widefield microscope can be modeled in intensity as a linear space-invariant system~\cite{Vonesch2006}, the adopted forward model in~\eqref{eq:observation_model} is still valid and we can, thus, employ the proposed framework for the restoration of the underlying biomedical images.

To evaluate the practical relevance of our approach, we provide deblurring experiments on two real images of fluorescence specimens. For each sample we acquired two image-stacks using a confocal microscope. This type of microscope can reject out-of-focus light using a small aperture in front of the detector and can thus avoid the blurring effect, but at the expense of more measurement noise. When the aperture is opened, the intensity of the incoming light is increased and the SNR is improved, but this time the measurements include interference from adjacent out-of-focus objects. In this case the final result is blurred and it is equivalent to an image acquired by a ``cheaper" widefield microscope. For more details on the image acquisition we refer to~\cite{Vonesch2008}.

\begin{figure}[!t]
\centering
 \includegraphics[scale=.35]{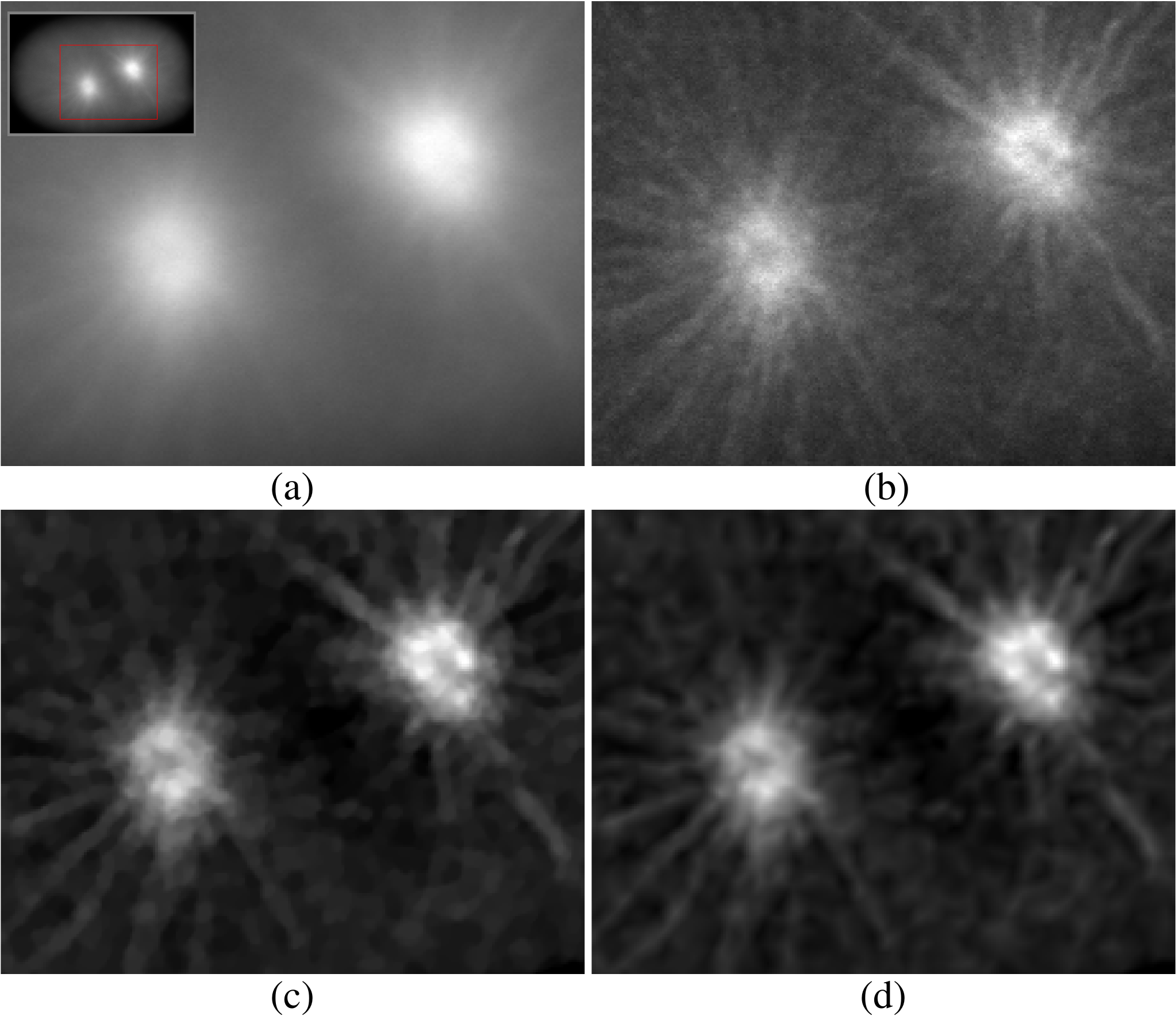}\vspace{-.2cm}
 \caption{Restoration results on a real fluorescent-cell image of size
$352\times 512$. Close-up of (a) widefield image, (b) reference confocal image,
(c) TV reconstruction, (d) $\mc{HS}_1$ reconstruction.
The details of this figure are better seen in the electronic version of
this paper by zooming on the screen.}
 \label{fig:bio1_results}
\end{figure}

\begin{figure}[!t]
\centering
 \includegraphics[scale=.45]{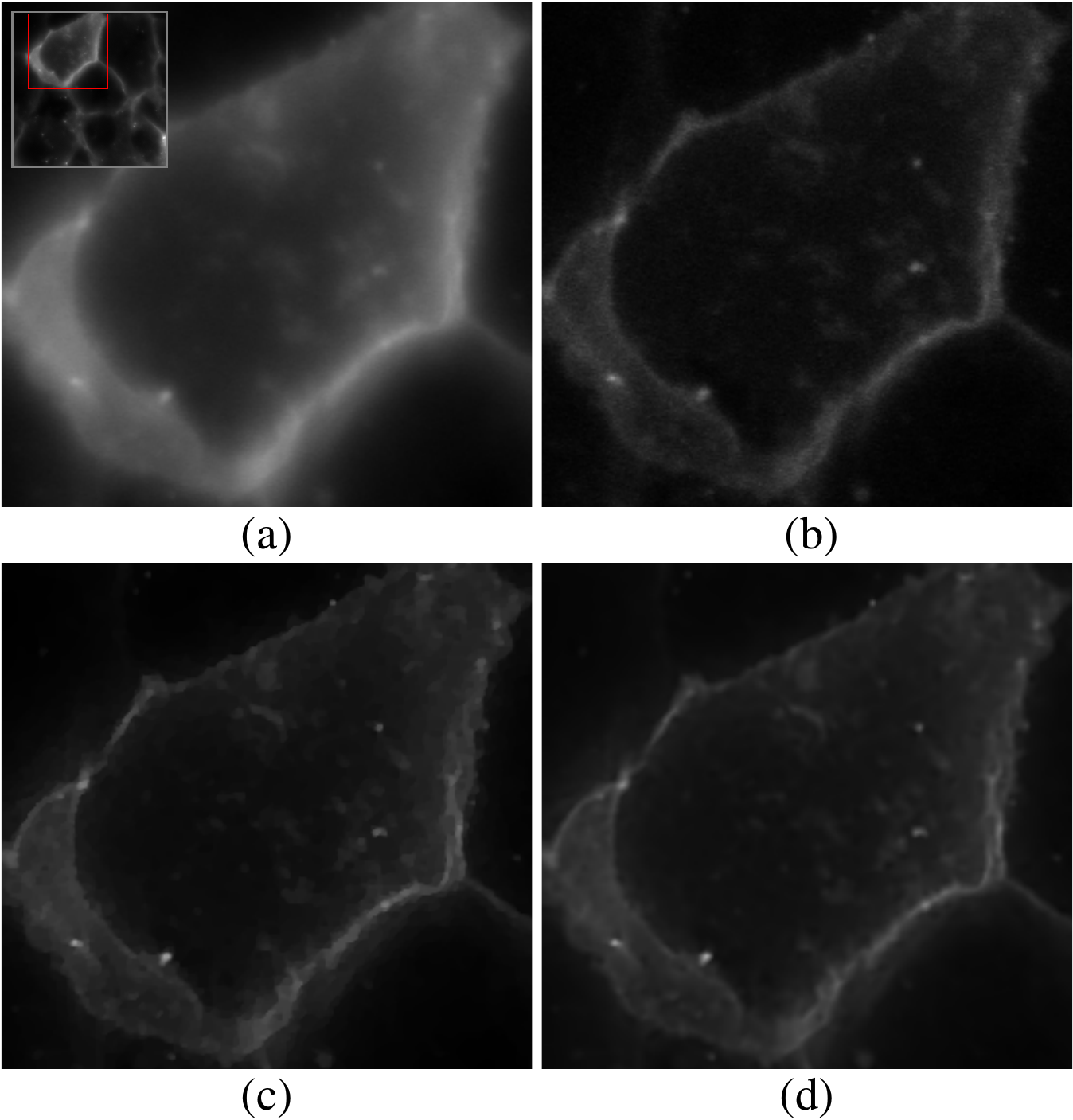} \vspace{-.2cm}
 \caption{Restoration results on a real fluorescent-cell image of size
$512\times 512$. Close-up of (a) widefield image, (b) reference confocal image,
(c) TV reconstruction, (d) $\mc{HS}_2$ reconstruction.
The details of this figure are better seen in the electronic version of
this paper by zooming on the screen.}
 \label{fig:bio2_results}
 \vspace{-.25cm}
\end{figure}

The reported results refer to the restoration of the second type of image-stacks, with the first ones serving as visual references to evaluate the quality of the reconstruction. The size of the image-stacks for the first specimen shown in Figs.~\ref{fig:bio1_results}(a) and \ref{fig:bio1_results}(b) are $352\times 512\times 96$ while the size of the image stacks for the second sample shown in Figs.~\ref{fig:bio2_results}(a) and \ref{fig:bio2_results}(b) are $512\times 512\times 16$. From each of the image-stacks we obtained a single image to work with, by computing the average intensity with respect to the z-axis. We did the same to obtain a 2D PSF out of a standard diffraction-limited 3D PSF model using the nominal optical parameters of the microscope (numerical aperture, wavelength, optical zoom)~\cite{Vonesch2006}.

In Figs.~\ref{fig:bio1_results}(c) and \ref{fig:bio1_results}(d) we present the restored images using TV and $\mc{HS}_1$ regularization, while in Figs.~\ref{fig:bio2_results}(c) and \ref{fig:bio2_results}(d) we provide the restored images using TV and $\mc{HS}_2$ regularization. From these two examples, if we compare the obtained results with the confocal acquisitions, we can verify that the Hessian-based solutions are quite successful in revealing the primary features of the specimens without introducing severe artifacts, as opposed to TV which oversmooths certain features and wipes out important details of the image structure.  Therefore, we conclude that our regularizers can do a better job, especially when one has to deal with images that consist mostly of ridges and filament-like structures, as is often the case in biomedical imaging.

\vspace{-.25cm}
\subsection{Sparse Image Reconstruction}
\label{sec:SparseRecon}

\begin{table}[!t]
 \centering
  \caption{PSNR comparisons on sparse image reconstruction from random samples for 4 ratios of observed pixels}
    \includegraphics[scale=.9]{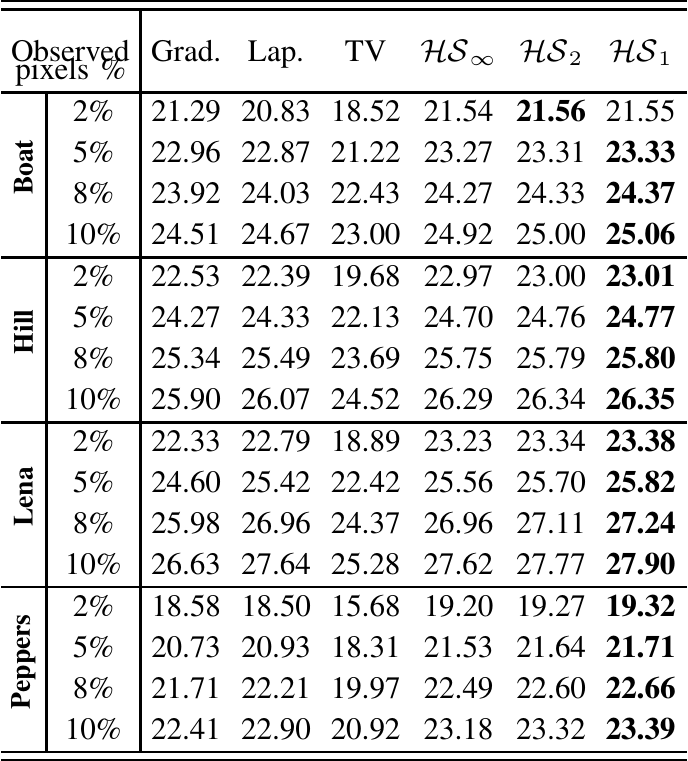}
 \label{tab:sparse_results}
 \vspace{-.5cm}
\end{table}

In sparse image reconstruction the observed image $\my$ is degraded by a masking operator  which randomly sets pixel values to zero. This operator corresponds to a diagonal matrix $\mA$ whose diagonal entries are randomly set to zero or one. We refer to this problem as sparse because in our experiments we consider masking operators that retain only $2\%, 5\%, 8\%$ and $10\%$ of the initial pixel values. Note that this problem can be considered as compressive sensing if we assume that the image is sparse in the spatial domain.  

The reported experiments are conducted on the gray-level images: Boat, Hill, Lena and Peppers. The masked images are then reconstructed using our regularizers as before plus TV and two quadratic regularizers based on the gradient and the Laplacian operators, respectively. 
In this setting we do not consider any presence of noise and thus, for all the regularizers under comparison, we use the same regularization parameter $\tau=10^{-4}$. The value of $\tau$ is chosen to be small to ensure that the results will be consistent, in the sense that the reconstruction methods will not alter the unmasked pixel values. However, due to the small value of the regularization parameter, we have observed that the convergence of the minimization task for all the regularizers can be slow and thus more than 100 iterations are required. To cope with this problem, for the non-quadratic regularizers, we apply a simple continuation scheme that significantly speeds up the convergence: We start with a large value for $\tau$ and then we gradually decrease it to reach the chosen value. We observe experimentally that following this strategy and using 200 MFISTA iterations (we still solve the corresponding denoising problems using 10 iterations) we can solve the problem to high accuracy. Regarding the two quadratic regularizers, we minimize their objective functions using the conjugate gradient method~\cite{Shewchuk1994} with a maximum of 2000 iterations. 

\begin{figure}[!t]
  \centering
    \includegraphics[scale=.4]{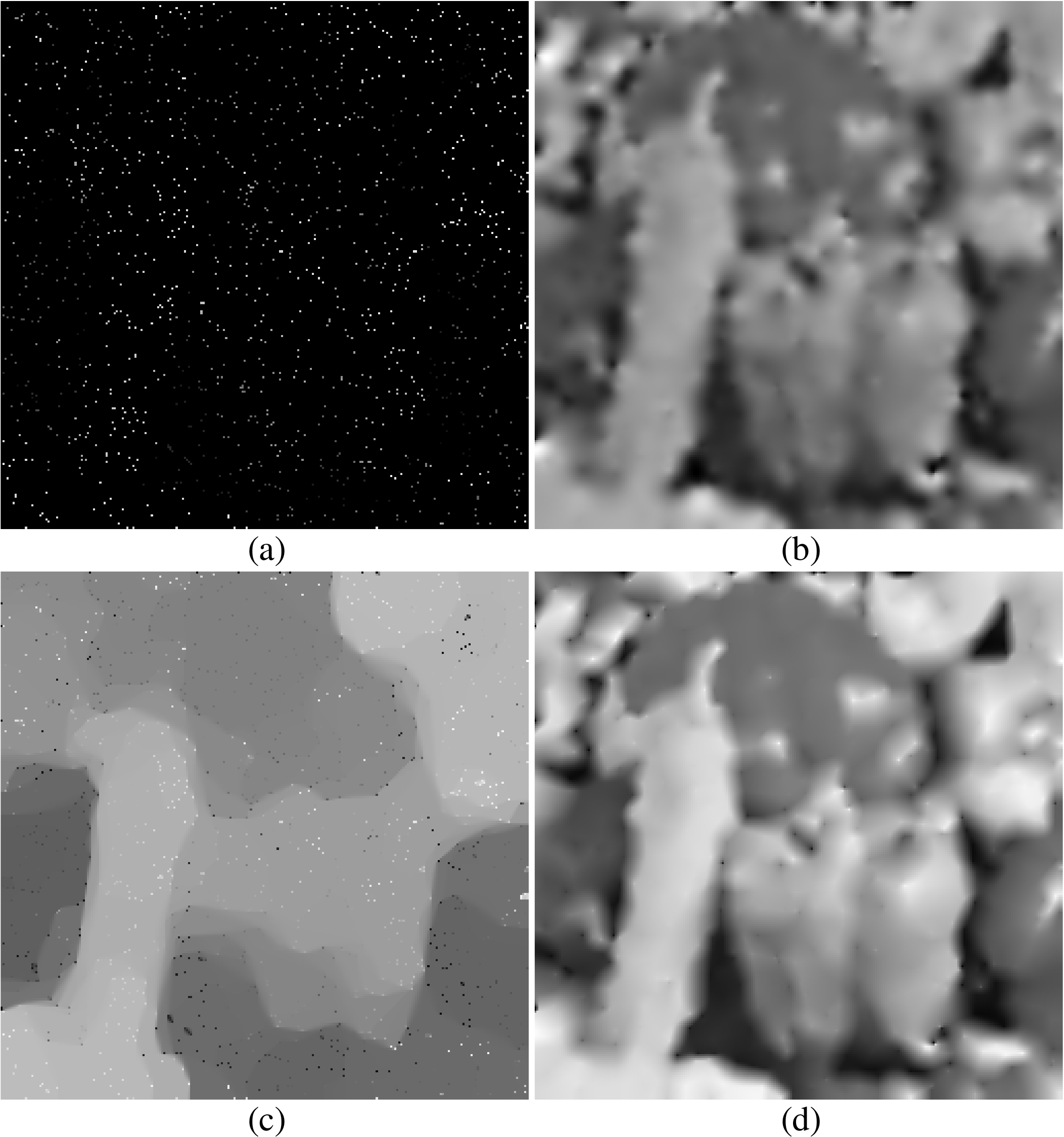}\vspace{-.2cm}
    \caption{Sparse reconstruction of the Peppers image from $2\%$ observed pixels. (a) Masked image,  (b) Laplacian-based quadratic result (PSNR = 18.50 dB),
 (c) TV result (PSNR = 15.68 dB),  and (d) $\mc{HS}_1$ result (PSNR = 19.32 dB)}
 \label{fig:peppers_random_reconstruction}
 \vspace{-.4cm}
\end{figure}

In Table~\ref{tab:sparse_results} we report the reconstruction results we obtained for all the employed regularization techniques. The quality of the reconstructions is measured in terms of PSNR. As we can observe from this table, TV does not perform well in this problem and its reconstructions fall far behind, even from the two quadratic regularization techniques. On the other hand, our Hessian-based regularizers behave much better and in all cases they lead to estimates that outperform the other methods. As in the image restoration case, the $\mc{HS}_1$ regularizer leads to the best reconstructions while the $\mc{HS}_2$ regularizer follows rather closely.  In Fig.~\ref{fig:peppers_random_reconstruction} we present a representative example of the reconstruction of the Peppers image from $2\%$ observed pixels. From this example it is clear that TV cannot produce an acceptable  result but instead leads to a piecewise constant solution that does not reveal any features of the image. On the other hand, both the quadratic and the proposed regularizer provide meaningful reconstructions with the latter achieving a better performance. 

\vspace{-.2cm}
\subsection{Image Interpolation and Image Zooming}
\label{sec:Zooming}
\begin{figure*}[!t]
  \centering
    \includegraphics[scale=.7]{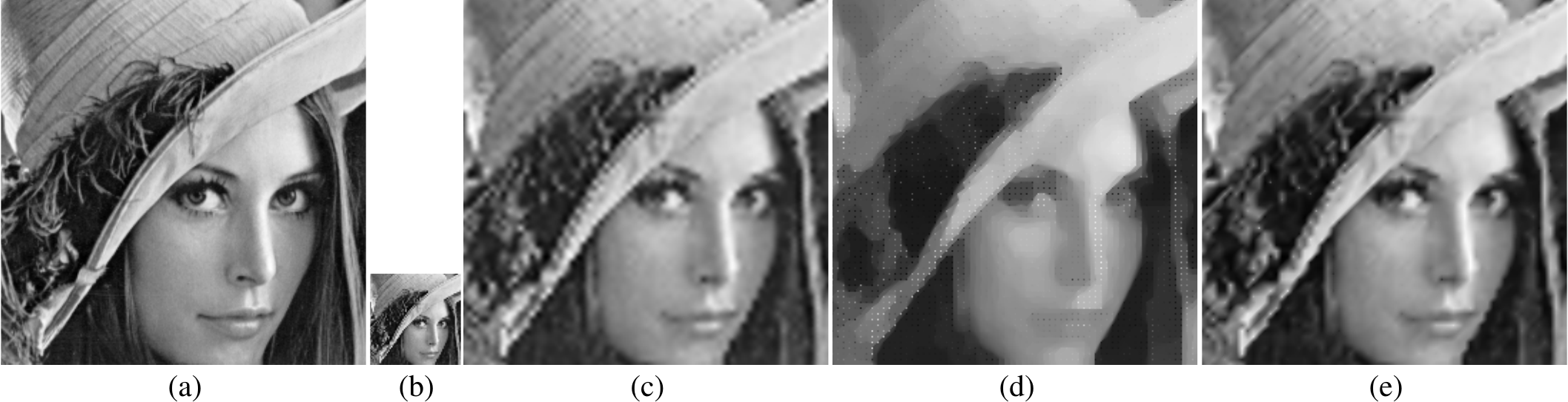}\vspace{-.2cm}
    \caption{Image interpolation. Close-up of (a) High-resolution image, (b) low-resolution image, (c) Laplacian-based quadratic result (PSNR=27.74 dB),  (d) TV result (PSNR=23.95 dB), and (e) $\mc{HS}_1$ result (PSNR=27.92 dB)}
 \label{fig:lena_interpolation}
\end{figure*}

\begin{figure*}[!t]
  \centering
    \includegraphics[scale=.7]{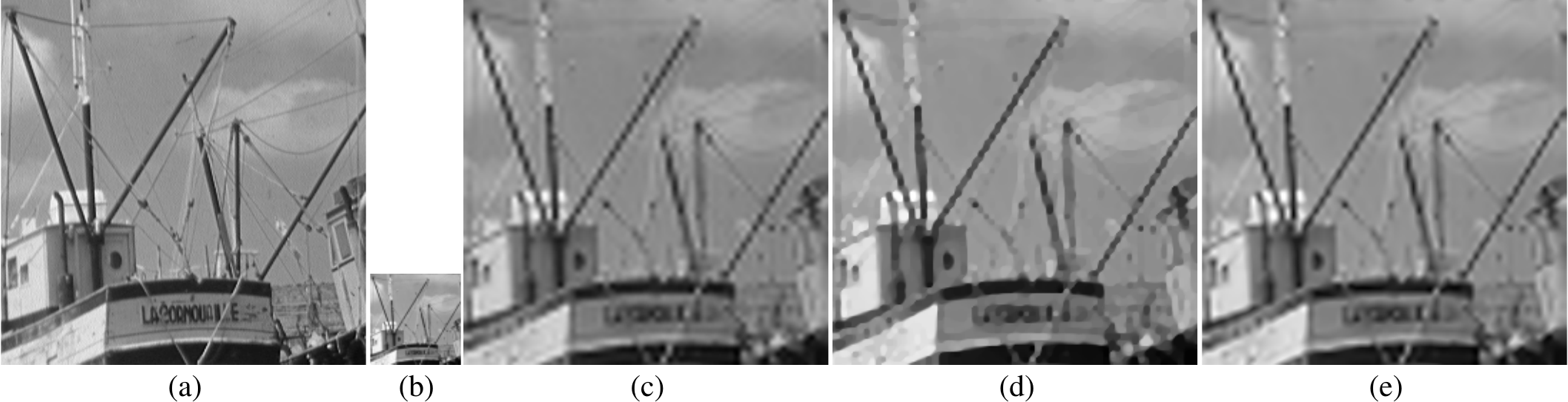}\vspace{-.2cm}
    \caption{Image zooming. Close-up of (a) High-resolution image, (b) low-resolution image,  (c) gradient-based quadratic result (PSNR = 25.95 dB),  (d) TV result (PSNR = 26.00 dB), and (e) $\mc{HS}_2$ result (PSNR = 26.08 dB)}
 \label{fig:boat_zooming}
 \vspace{-.4cm}
\end{figure*}

Image interpolation and image zooming fall into the same class of linear inverse problems. As in the sparse image reconstruction case,  the degradation is due to a masking operator that zeros out some of the image pixel values. However, in these two cases the masking operator corresponds to subsampling and is highly structured, as opposed to the random masking operator. The difference between the two considered forward models is that image interpolation involves only the subsampling operator and therefore results in observed images that suffer from aliasing, while image zooming  involves additionally  an antialiasing operator that is applied to the underlying image before the subsampling takes place. In the last case, the matrix $\mA$ can be expressed as $\mA=\m{S}\m{F}$ where $\m{F}$ corresponds to the filtering operation and $\m{S}$ to subsampling.  Once more, we do not consider any presence of noise and we thus use the same regularization parameter and minimization strategy as above. The experiments we present are conducted on the same four images as in Section~\ref{sec:SparseRecon}, using the same regularizers for a downsampling factor of 4. Finally, regarding the antialiasing filter we use a Gaussian kernel of support $9\times 9$ and standard deviation $\sigma_b=1.4$.

In Table~\ref{tab:zooming_results} we report the obtained results and we evaluate the quality of the estimates in terms of PSNR.  Regarding the interpolation problem we observe that TV, similarly to the sparse image reconstruction case, does not perform well and produces the worst scores. However, its performance gets significantly better in the image zooming case where an antialising filtering is applied. This is an indication that TV cannot perform at a satisfactory level when the operator acting on the image does not involve a mixing effect.  On the other hand, the performance of the proposed regularizers is more robust to the nature of the degradation operator, and they lead to the best reconstructions. To have a visual performance assessment, we present in Figs.~\ref{fig:lena_interpolation} and \ref{fig:boat_zooming} interpolation and zooming results on the Lena and Boat image, respectively. These results confirm our previous conclusions about the performance of TV and demonstrate the superiority of the Hessian-based regularizers over the other regularizers.

\begin{table}[]
 \centering
  \caption{PSNR comparisons on image interpolation and image zooming for a $4\times$ downsampling factor}
    \includegraphics[scale=.85]{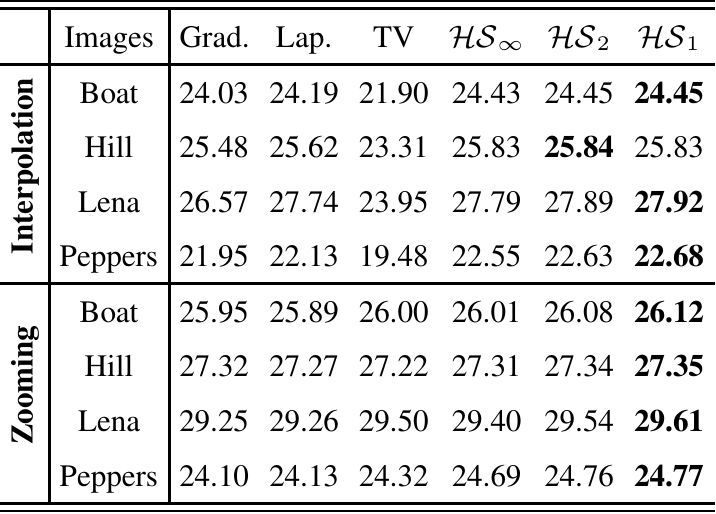}
 \label{tab:zooming_results}
 \vspace{-.2cm}
\end{table}

\section{Conclusion}
\label{sec:conclusions}
In this paper we introduced a new family of convex non-quadratic regularizers that can potentially lead to improved results in inverse imaging problems. These regularizers incorporate second-order information of the image and depend on the Schatten norms of the Hessian. We further designed an efficient and highly parrallelizable projected gradient algorithm for minimizing the corresponding objective functions. We also presented a new result that relates vector projections onto $\ell_q$ norm balls and matrix projections onto Schatten norm balls. This enabled us to design a matrix-projection method, which is a fundamental ingredient of our optimization algorithm.

The performance and practical relevance of the proposed regularization scheme was assessed for several linear inverse imaging problems, through comparisons on simulated and real experiments with various competing methods. The results we obtained are promising and competitive both in terms of SNR improvement and visual quality.

\appendices
\section{}
\label{sec:appendix_A}
\subsection{Proof of Theorem~\ref{thm:grad_norm}}
\label{proof:grad_norm}
By taking the domain $\W$ to be a disk, the rotation invariance of $\mc{R}\br{f}$ implies that
\beqn
\mc{R}\br{f\br{\m{R}_\theta\cdot}}=\mc{R}\br{f}
\label{property:Rot_invar_grad}
\eeqn
where $\m{R}_\theta$ is a rotation matrix. In particular, \eqref{property:Rot_invar_grad} must hold for all functions, including those of the form:
$f\br{\m{r}}=\a r_1+\b r_2$, with $\mr\in\R^2$ and $\a,\b\in \R$. Their gradient is constant and equal to $\grad f\br{\mr}=\pr{\bsmtx \a\\ \b\esmtx}=\mx$.
Now, using $f$ as defined above, we write the l.h.s of~\eqref{property:Rot_invar_grad} as
\begin{align}
\mc{R}\br{f\br{\m{R}_\theta\cdot}}&=\int_{{\W}}\Phi\br{\grad\cbr{f\br{\m{R}_\theta\cdot}}\br{\mr}}\mbox{d}\mr\nonumber\\
&=\int_{{\W}}\Phi\br{\m{R}_{\theta}^\transpose\grad f\br{\m{R}_\theta\mr}}\mbox{d}\mr
=\int_{{\W}}\Phi\br{\m{R}_{\theta}^\transpose\mx}\mbox{d}\mr\nonumber\\
&=\int_{{\W}}\Phi\br{\abs{\mx}\cdot\m{u}_{\t'}}\mbox{d}\mr\,,
\label{eq:rot_inv_grad_eq}
\end{align}
where $\m{u}_{\t'}=\bbmtx\sin\br{\t'}& \sin\br{\t'+\frac{\pi}{2}}\ebmtx^\transpose$ and $\t'=\t+\sgn{\a}\arccos\br{\frac{\b}{\sqrt{\a^2+\b^2}}}$. Setting $\t'=\frac{\pi}{2}$ in \eqref{eq:rot_inv_grad_eq} and combining the result with~Property~\eqref{property:Rot_invar_grad}, we immediately get that
\begin{align}
\Phi\br{\mx}=\Phi\br{\abs{\mx}}, \forall \mx\in\R^2\,.
\label{eq:Phi_grad}
\end{align}
The scaling invariance of $\mc{R}$ can be restated as 
\begin{align}
\mc{R}_a\br{f\br{a\cdot}}&=a^{\mu}\mc{R}\br{f}\nonumber\\
\int_{\W/a}\Phi\br{\grad\cbr{f\br{a\cdot}}\br{\mr}}\mb{d}\mr&=a^\mu\int_{\W}\Phi\br{\grad f\br{\mr}}\mb{d}\mr
\label{property:Scale_invar_grad}
\end{align}
for some $a>0$, and an exponent $\mu\in\R$. This property must hold for all functions, including those of the form:
$f\br{\m{r}}=\a r_1$, with $\mr\in\R^2$ and $\a\in\R$. The magnitude of their gradient is constant and equal to 
$\abs{\grad f\br{\mr}}=\abs{\a}$.
Now, using $f$ as defined above and the result of~\eqref{eq:Phi_grad}, we write~\eqref{property:Scale_invar_grad} as
\begin{align}
\int_{\W/a}\Phi\br{a\abs{\a}}\mb{d}\mr=a^\mu\int_{\W}\Phi\br{\abs{\a}}\mb{d}\mr\,.
\label{eq:homogeneity_grad}
\end{align}
Therefore, we directly have that 
\beqn
\Phi\br{a\abs{\a}}=a^\nu\Phi\br{\abs{\a}}, \forall \a\in\R,
\label{eq:scale_invar_Phi}
\eeqn 
with $\nu=\mu+2$. Now, we define the function
\beqn
\Phi_0(\a)=\frac{\Phi(\abs{\a})}{\abs{\a}^\nu}, \forall \a\in\R
\eeqn
which is homogeneous of degree 0. This implies that $\Phi_0(\a)=c$, with $c$ an arbitrary constant. Therefore, the potential functions $\Phi$  satisfying~\eqref{eq:scale_invar_Phi} are necessarily of the form: $\Phi\br{\cdot}=c\abs{\cdot}^\nu$. 

The inverse statement of the theorem can be verified by substitution, using the property 
\beqn
\grad\cbr{f\br{a\cdot}}\br{\mr}=a\grad\cbr{f}\br{a\mr},\,\forall f.
\eeqn

\subsection{Proof of Theorem~\ref{thm:Hessian_norm}}
\label{proof:Hessian_norm}
By taking the domain $\W$ to be a disk, the rotation invariance of $\mc{R}\br{f}$, as defined in~\eqref{property:Rot_invar_grad},
must hold for all functions, including those of the form:
$f\br{\m{r}}=\frac{\a}{2} r_1^2+\frac{\b}{2} r_2^2+\gamma r_1 r_2$, 
with $\mr\in\R^2$ and $\a,\b,\gamma\in \R$. Their Hessian is constant and equal to $\mc{H}f\br{\mr}=\pr{\bsmtx \a & \gamma \\ \gamma& \b \esmtx}=\mA$. Now, using $f$ as defined above, we write the l.h.s of~\eqref{property:Rot_invar_grad} as
\begin{align}
\mc{R}\br{f\br{\m{R}_\theta\cdot}}&=\int_{\W}\Phi\br{\mc{H}\cbr{f\br{\m{R}_\theta\cdot}}\br{\mr}}\mbox{d}\mr\nonumber\\
&=\int_{\W}\Phi\br{\m{R}_{\theta}^\transpose\mc{H}f\br{\m{R}_\theta\mr}\m{R}_{\theta}}\mbox{d}\mr\nonumber\\
&=\int_{\W}\Phi\br{\m{R}_{\theta}^\transpose\mA\m{R}_{\theta}}\mbox{d}\mr\,.
\end{align}
According to the spectral decomposition theorem, $\mA$ being symmetric has an eigenvalue decomposition. 
This implies that there exists a rotation $\theta'$ such that $\m{R}_{\theta'}^\transpose\mA\m{R}_{\theta'}=\m{\Lg}$, where $\m{\Lg}$ is a diagonal matrix consisting of the eigenvalues of $\mA$. These are denoted as $\m{\lg}_k$, where $k=1,2$.  Based on this observation and combining it with Property~\eqref{property:Rot_invar_grad}, we immediately get that
\begin{align}
\Phi\br{\mA}=\Phi\br{\m{\lg}_1,\m{\lg}_2},\forall \mA\in \S^2,
\label{eq:Phi_eigen}
\end{align}
which implies that $\Phi$ should be a function of the Hessian eigenvalues. 

The scaling invariance of $\mc{R}$ can be restated as 
\begin{align}
\mc{R}_a\br{f\br{a\cdot}}&=a^{\mu}\mc{R}\br{f}\nonumber\\
\int_{\W/a}\Phi\br{\mc{H}\cbr{f\br{a\cdot}}\br{\mr}}\mb{d}\mr&=a^\mu\int_{\W}\Phi\br{\mc{H}f\br{\mr}}\mb{d}\mr
\label{property:Scale_invar_eigen}
\end{align}
for some $a>0$, and an exponent $\mu\in\R$. This property must hold for all functions, including those of the form:
$f\br{\m{r}}=\frac{\a}{2} r_1^2+\frac{\beta}{2} r_2^2$, with $\mr\in\R^2$ and $\a,\b\in\R$. Their Hessian is constant and equal to $\mc{H}f\br{\mr}=\pr{\bsmtx\a & 0\\ 0& \b \esmtx}$.
Now, using $f$ as defined above and the result of~\eqref{eq:Phi_eigen}, we write~\eqref{property:Scale_invar_eigen} as
\begin{align}
\int_{\W/a}\Phi\br{a^2\a,a^2\b}\mb{d}\mr=a^\mu\int_{\W}\Phi\br{\a,\b}\mb{d}\mr\,.
\label{eq:homogeneity_eigen}
\end{align}
Therefore, we have that 
\beqn
\Phi\br{a^2\a,a^2\b}=a^{\mu+2}\Phi\br{\a,\b},\, \forall \a,\b\in\R\,.
\label{eq:scale_invar_Phi_eigen}
\eeqn 
Now, we define the function 
\beqn
\Phi_0(\mx)=\frac{\Phi\br{\mx}}{\norm{\mx}{p}^{\nu}}\,,\forall \mx\in\R^2,
\eeqn
where $p\ge 1$ and $\nu=\frac{\mu+2}{2}$. $\Phi_0$ is homogeneous of degree 0 and thus $\Phi_0\br{\mx}=\Phi_0\br{\mx/\norm{\mx}{p}}$. Therefore, the potential functions $\Phi$ that satisfy \eqref{eq:scale_invar_Phi_eigen}, are necessarily of the form: $\Phi\br{\mx}=\Phi_0\br{\mx/\norm{\mx}{p}}\norm{\mx}{p}^{\nu}$. Finally, since $\mx$ represents the vector of the eigenvalues of the Hessian, its $\ell_p$ norm corresponds to the $\mc{S}_p$ norm of the Hessian itself. 

The inverse statement of the theorem can be verified by substitution, using the property 
\beqn
\mc{H}\cbr{f\br{a\cdot}}\br{\mr}=a^2\mc{H}\cbr{f}\br{a\mr},\,\forall f.
\eeqn

\section{}
\label{sec:appendix_B}
\subsection{Adjoint of the Disrcete Hessian Operator}
\label{sec:Adjoint_Hessian}
To find the adjoint of the discrete Hessian operator, we exploit the relation of the inner products of the spaces $\R^N$ and $\mc{X}$ in~\eqref{eq:inner_product_relation}. Using~\eqref{eq:inner_product_of_X},  we can equivalently write~\eqref{eq:inner_product_relation} as 
\beqn
\sum_{n=1}^N\tr{\pr{\dH\mx}_n^\transpose\mY_n}=\sum_{n=1}^N x_n\pr{\dH^*\mY}_n\,.
\label{eq:inner_product_relation_discrete}
\eeqn
We then expand the l.h.s of \eqref{eq:inner_product_relation_discrete},  to obtain
\begin{align}
&\sum_{n=1}^N\tr{\pr{\dH\mx}_n^\transpose\mY_n}=\sum_{n=1}^N\Big(\pr{\dH\mx}_n^{\br{1,1}}\mY_n^{\br{1,1}}+\nonumber\\&\quad\quad\quad\pr{\dH\mx}_n^{\br{1,2}}\br{\mY_n^{\br{1,2}}+\mY_n^{\br{2,1}}}+\pr{\dH\mx}_n^{\br{2,2}}\mY_n^{\br{2,2}}\Big)\nonumber\\
=&\sum_{n=1}^N\Big(\pr{\Dxx\mx}_n\mY_n^{\br{1,1}}+\pr{\Dxy\mx}_n\br{\mY_n^{\br{1,2}}+\mY_n^{\br{2,1}}}
\nonumber\\&\quad\quad\quad+\pr{\Dyy\mx}_n\mY_n^{\br{2,2}}\Big)\nonumber\\
=&\sum_{n=1}^N x_n\Big(\pr{\Dxx^*\mY^{\br{1,1}}}_n+\pr{\Dxy^*\br{\mY^{\br{1,2}}+\mY^{\br{2,1}}}}_n
\nonumber\\ &\quad\quad\quad+\pr{\Dyy^*\mY^{\br{2,2}}}_n\Big)\,.
\label{eq:expansion_of_inner_product_of_X}
\end{align}
Note that $\mY^{\br{i,j}}$ corresponds to the vector that is composed of the $\br{i,j}$ entries of all $\mY_n\in\R^{2\times 2}$ matrices. Now,  by comparing the r.h.s of~\eqref{eq:inner_product_relation_discrete} to the r.h.s expansion of~\eqref{eq:expansion_of_inner_product_of_X}, it is straightforward to verify that the adjoint of the discrete Hessian operator is indeed computed according to~\eqref{eq:Hessian_adjoin_Op}.

\subsection{Proof of Proposition~\ref{thm:matrix_projection}}
\label{proof:matrix_projection}
By definition, the orthogonal projection of a matrix $\mY$ onto the set $\mc{B}_{\mc{S}_q}$ is given by
\beqn
\mc{P}_{\mc{B}_{\mc{S}_q}}\br{\mY}=\argmin_{\snorm{\mX}{q}\le \rho}
\norm{\mX-{\mY}}{F}^2\,.
\label{eq:Proj_def}
\eeqn
Since all Schatten norms are unitarily invariant, we equivalently have
\beqn
\mc{P}_{\mc{B}_{\mc{S}_q}}\br{\mY}=\argmin_{\snorm{\mU^\hermitian\mX\mV}{q}\le \rho}
\norm{\mU^\hermitian\mX\mV-\mU^\hermitian{\mY}\mV}{F}^2\,.
\label{eq:Proj_def_equiv}
\eeqn
Let us now consider the matrix $\mZ=\mU^\hermitian\mX\mV$ that is associated with the solution of \eqref{eq:Proj_def}.  If we substitute $\mZ$ in \eqref{eq:Proj_def_equiv}, then we end up with the following constrained minimization problem
\beqn
\mc{P}_{\mc{B}_{\mc{S}_q}}\br{\mS}=\argmin_{\snorm{\mZ}{q}\le \rho}
\norm{\mZ-\mS}{F}^2\,,
\label{eq:Proj_def_diag}
\eeqn
which corresponds to the projection of the diagonal matrix $\mS$ onto the set $\mc{B}_{\mc{S}_q}$. Now, if $\mc{P}_{\mc{B}_{\mc{S}_q}}\br{\mS}=\hat{\mZ}$, we have
\begin{align}
\norm{\hat{\mZ}-\mS}{F}^2&=\norm{\hat{\mZ}}{F}^2+\norm{\mS}{F}^2-2\operatorname{Re}\br{\tr{\hat{\mZ}^\hermitian{\mS}}}\nonumber\\
&\ge\norm{\hat{\mS}}{F}^2+\norm{\mS}{F}^2-2\tr{\hat{\mS}^\transpose\mS}=\norm{\hat{\mS}-\mS}{F}^2\,,
\label{eq:Fro_distance}
\end{align}
where the inequality stems from von Neumann's trace theorem~\cite{Mirsky1975}, and $\hat{\mS}$ is a diagonal matrix with the singular values of $\hat{\mZ}$. In addition, it holds that 
\beqn
\snorm{\hat{\mS}}{q}=\snorm{\hat{\mZ}}{q} \le \rho\,.
\label{eq:schatten_constraint}
\eeqn
Equations \eqref{eq:Fro_distance} and \eqref{eq:schatten_constraint} immediately imply that the projection of $\mS$ equals to $\hat{\mZ}=\hat{\mS}$, i.e., a positive semidefinite diagonal matrix.  We can then perform this operation by projecting the vector, formed by the main diagonal of $\mS$, onto the convex set $\mc{B}_q$, and then by transforming the projected vector back to a diagonal matrix. Using this fact and the relations between the optimal solution of \eqref{eq:Proj_def} and \eqref{eq:Proj_def_diag}, we finally express the projection of the matrix $\mY$ onto $\mc{B}_{\mc{S}_q}$ as
\begin{align}
\mc{P}_{\mc{B}_{\mc{S}_q}}\br{\mY}=\mU\diag{\mc{P}_{\mc{B}_q}\br{\ms\br{\mY}}}\mV^\hermitian\,.
\label{eq:Mat_Proj}
\end{align}

\subsection{Proof of Lemma~\ref{lemma:lemma_dual_norm}}
\label{sec:Proof_lemma_dual_norm}
First, we present a matrix inequality that involves the Schatten norms and it will be subsequently used for the proof of the lemma. Let $\mX$, $\mY\in\C^{n_1\times n_2}$. Then, the inner product of these two matrices satisfies the following inequality
\begin{align}
\ip{\mX}{\mY}_{\C^{n_1\times n_2}}&=\operatorname{Re}\br{\tr{\mY^\hermitian\mX}}\leq\ip{\ms\br{\mX}}{\ms\br{\mY}}_2\nonumber\\
&\leq\norm{\ms\br{\mX}}{q}\norm{\ms\br{\mY}}{p}=\snorm{\mX}{q}\snorm{\mY}{p}\,.
\label{eq:Schatten_matrix_inequality}
\end{align}
The first inequality is due to  von Neumann's trace theorem \cite{Mirsky1975}, while the second one due to H\"{o}lder's inequality. The last equality holds true from the definition of Schatten norms. 

By definition, the dual norm of~\eqref{eq:mixed_l1_Sp} is given by~\cite{Rockafellar1970}:
\beqn
\norm{\m{\W}}{D}=\max_{\norm{\m{\P}}{1,p}\leq 1}\ip{\m{\W}}{\m{\P}}_{\mc{X}}\,,
\label{eq:dual_norm_definition}
\eeqn
where $\mc{X}$, instead of $\R^{N\times 2\times 2}$ that is used throughout  the paper, here is assumed to be the more general linear space $\mc{X}=\C^{N\times n_1\times n_2}$.
We consider the inequality 
\begin{align}
\ip{\m{\W}}{\m{\P}}_{\mc{X}}=\sum_{n=1}^N\operatorname{Re}\br{\tr{\mP_n^\hermitian\mW_n}}\leq\sum_{n=1}^N\snorm{\mW_n}{q}\snorm{\mP_n}{p}\,,
\label{eq:inequality_inner_product_of_X_a}
\end{align}
which immediatelly follows from inequality~\eqref{eq:Schatten_matrix_inequality}. Now, by introducing the vectors $\bm{\w}=\br{\snorm{\mW_1}{q}, \snorm{\mW_2}{q},\ldots,\snorm{\mW_N}{q}}$ and $\bm{\p}=\br{\snorm{\mP_1}{p}, \snorm{\mP_2}{p},\ldots,\snorm{\mP_N}{p}}$, and applying once again H\"{o}lder's inequality, we get
\begin{align}
\sum_{n=1}^N\snorm{\mW_n}{q}\snorm{\mP_n}{p}&=\ip{\bm{\w}}{\bm{\p}}_2\leq\norm{\bm{\w}}{\infty}\norm{\bm{\p}}{1}\nonumber\\&=\norm{\mW}{\infty,q}\norm{\mP}{1,p}\,.
\label{eq:inequality_inner_product_of_X_b}
\end{align}
From the definition of the dual norm~\eqref{eq:dual_norm_definition} and the inequalities~\eqref{eq:inequality_inner_product_of_X_a} and~\eqref{eq:inequality_inner_product_of_X_b}  we conclude that $\norm{\m{\W}}{D}\leq\norm{\m{\W}}{\infty,q}$. To prove that $\norm{\m{\W}}{D}=\norm{\m{\W}}{\infty,q}$, we next show that for each $\mW$ we can find a $\mP$ satisfying $\norm{\mP}{1,p}=1$, and for which $\ip{\m{\W}}{\m{\P}}_{\mc{X}}=\norm{\m{\W}}{\infty,q}$. To that end, let $k$ be any index in the set $\cbr{\argmax_{1\le n\le N}\snorm{\mW_n}{q}}$ and  $\mW_k=\m{U}_k\mS_k\m{V}^\hermitian_k$ be the singular value decomposition of $\mW_k$. Then, we set $\mP_n=\m{O}$ for all $n$ except for $n=k$ for which we have
\beqn
\mP_k=\m{U}_k\m{E}\m{V}^\hermitian_k\,,
\eeqn
where 
\beqn
\m{E}^{\br{i,j}}=\frac{\br{\mS^{\br{i,j}}_k}^{q-1}}{\snorm{\mW_k}{q}^{q-1}}\,,
\eeqn
and $\mS^{\br{i,j}}_k$ corresponds to the $\br{i,j}$-th entry of the matrix $\mS_k\in\D^{n_1\times n_2}$. Now, we have that 
\begin{align}
\ip{\mW}{\mP}_{\mc{X}}&=\sum_{n=1}^N\operatorname{Re}\br{\tr{\mP_n^\hermitian\mW_n}}=\operatorname{Re}\br{\tr{\mP_k^\hermitian\mW_k}}\nonumber\\
&=\tr{\m{E}^\hermitian\mS_k}=\frac{\suml_{i=1}^{\min\br{n_1,n_2}}\br{\mS^{\br{i,i}}_k}^q}{\snorm{\mW_k}{q}^{q-1}}\nonumber\\
&=\snorm{\mW_k}{q}=\norm{\mW}{\infty,q}\,.
\end{align}
Furhermore, for the mixed norm $\norm{\mP}{1,p}=\snorm{\mP_k}{p}$ it holds
\begin{align}
\norm{\mP}{1,p}&=\br{\suml_{i=1}^{\min\br{n_1,n_2}}\br{\m{E}^{\br{i,i}}}^p}^{1/p}
\nonumber\\
&=\frac{\br{\suml_{i=1}^{\min\br{n_1,n_2}}\br{\mS^{\br{i,i}}_k}^{q}}^{\frac{q-1}{q}}}{\snorm{\mW_k}{q}^{q-1}}=\frac{\snorm{\mW_k}{q}^{q-1}}{\snorm{\mW_k}{q}^{q-1}}=1\,,
\end{align}
which completes the proof of the lemma.

\subsection{Proof of Proposition~\ref{prop:Lipschitz_constant}}
\label{proof:Lipschitz_constant}
For any pair of variables $\mW\,,\mP\,\in \mc{X}$ we have
\begin{align}
\norm{\grad{s}\br{\mW}-\grad{s}\br{\mP}}{\mc{X}}&=\norm{\tau\dH\br{\mc{V}\br{\mW}-
\mc{V}\br{\mP}}}{\mc{X}}\nonumber\\
&\le\tau\norm{\dH}{}\norm{\mc{V}\br{\mW}-
\mc{V}\br{\mP}}{2}\nonumber\\
&\le\tau\norm{\dH}{}\norm{\tau\dH^*\br{\mW-\mP}}{2}\nonumber\\
&\le\tau^2\norm{\dH}{}\norm{\dH^*}{}\norm{\mW-\mP}{\mc{X}}\nonumber\\
&=\tau^2\norm{\dH}{}^2\norm{\mW-\mP}{\mc{X}},
\end{align}
where $\mc{V}\br{\mW}=\mc{P}_{\mc{C}}\br{\mz-\tau\dH^*\mW}$.
Note that, the first and third inequalities follow from the relation between the norms, defined in the spaces $\mc{X}$ and $\R^N$, and the induced operator norm, i.e., $\norm{\dH\mx}{\mc{X}}\le \norm{\dH}{}\norm{\mx}{2}$, while the second one holds because the projection operator $\mc{P}_\mc{C}$ onto the convex set $\mc{C}\subseteq\R^N$, is firmly nonexpansive~\cite[Proposition 4.8]{Bauschke2011}. This means that
\beqn
\norm{\mc{P}_\mc{C}\br{\mx}-\mc{P}_\mc{C}\br{\my}}{2}\le\norm{\mx-\my}{2}\,\forall\,{\mx,\my}\in\R^N\,.
\eeqn
To compute an upper bound of $\norm{\dH}{}$, we exploit that $\norm{\dH}{}^2=\norm{\dH^*\dH}{}$~\cite{Hutson2005} (a general property of bounded linear operators), and we get
\begin{align}
\norm{\dH^*\dH\mx}{2}&=\norm{\br{\Dxx^*\Dxx+2\Dxy^*\Dxy+\Dyy^*\Dyy}\mx}{2}\nonumber\\
&\le\br{\norm{\Dxx}{}^2+2\norm{\Dxy}{}^2+\norm{\Dyy}{}^2}\norm{\mx}{2}\,.
\end{align}
Now, using the definitions of the second-order differential operators in~\eqref{eq:differential_operators}, it is easy to show that each of $\norm{\Dxx}{}$, $\norm{\Dxy}{}$ and  $\norm{\Dyy}{}$ is smaller than or equal to 4. This immediately implies that
$\norm{\dH}{}\le 8$ and, hence, an upper bound of the Lipschitz constant of $\grad s\br{\mW}$ will be $L\br{s}\le \tau^2\norm{\dH}{}^2\le 64\tau^2$.

\section*{Acknowledgment}
\label{sec:acknowledgments}
The authors would like to acknowledge Jean-Charles Baritaux and Pouya Tafti for fruitful discussions. They would also like to thank the anonymous reviewers and the associate editor for their useful comments and suggestions.

\ifCLASSOPTIONcaptionsoff
  \newpage
\fi


\end{document}